\pdfoutput=1
\documentclass[11pt]{article}
\usepackage[T1]{fontenc}
\usepackage{lmodern,amsmath,amssymb,amsthm,mathtools,booktabs,array,graphicx,microtype}
\usepackage[a4paper,margin=26mm]{geometry}
\usepackage[hidelinks]{hyperref}
\allowdisplaybreaks
\newcommand{\R}{\mathbb R}

\newtheorem{theorem}{Theorem}[section]
\newtheorem{lemma}[theorem]{Lemma}
\newtheorem{proposition}[theorem]{Proposition}
\newtheorem{corollary}[theorem]{Corollary}
\theoremstyle{remark}

\numberwithin{equation}{section}
\hypersetup{pdftitle={The linear stability diagram of the elliptic Lagrange triangle},pdfauthor={Xijun Hu and Li Wu},pdflang=en}
\author{Xijun Hu\thanks{Partially supported by the National Natural Science
Foundation of China (Grant No.~12521001) and the Taishan Scholars Climbing
Program of Shandong Province (Grant No.~TSPD20240802).}
\and Li Wu\thanks{Partially supported by the National Natural Science
Foundation of China (Grant No.~12171281).}}
\date{}
\title{The linear stability diagram of the elliptic Lagrange triangle}
\begin{document}
\maketitle
\begin{abstract}
We determine the essential Floquet spectrum throughout the open stability
diagram of the elliptic Lagrange triangle together with the exact Jordan
canonical form on every interior spectral boundary. On the upper antiperiodic
degeneracy curve, a unique quadratic contact with an explicit Riccati arc
locates the eccentricity at which the boundary of hyperbolicity switches
branches, while along that arc the trace of the return on its invariant
Lagrangian graph decreases strictly with eccentricity. Globally on $[0,1)$, the
lower antiperiodic curve is strictly decreasing and the upper curve admits a
unique nondegenerate maximum; toward collision their mass parameters vanish at
orders $5/2$ and $1/4$ in $1-e$, respectively, with convergent power--log
expansions.
\end{abstract}

\section{Introduction and main results}
We determine the complete linear stability diagram of the elliptic equilateral
three-body solutions, identifying the spectrum of the essential period return
in every open region of the mass--eccentricity rectangle together with the
exact Jordan normal form along every interior spectral boundary. A central role
is played by the antiperiodic degeneracy locus, consisting of those parameters
at which the essential shape equation admits a nontrivial perturbation whose
position and velocity both reverse sign after one revolution; globally, its
two analytic branches govern both the loss of spectral stability and the
branch switch along the boundary of hyperbolicity.

For three positive masses define
\[
 \beta=\frac{27(m_1m_2+m_2m_3+m_3m_1)}{(m_1+m_2+m_3)^2},
 \qquad \beta=4\kappa(3-\kappa),\quad 0\le\kappa\le 3/2 ,
\]
including the zero-mass-parameter endpoint by continuity. After removing
translations and the homographic degrees of freedom, the two real position
coordinates satisfy, in forward true anomaly,
\begin{equation}\label{ii:model}
 \begin{gathered}
 \mathcal A_{\kappa,e}u=-u''-2Ju'+q^{-1}D_\kappa u=0,\\
 J=\begin{pmatrix}0&-1\\1&0\end{pmatrix},\quad
 D_\kappa=\operatorname{diag}(3-\kappa,\kappa),\qquad q(t)=1+e\cos t .
 \end{gathered}
\end{equation}
Danby studied linear stability first for the elliptic restricted problem \cite{DanR} and then for three finite masses \cite{DanG}. Roberts obtained the four-dimensional essential system, continued its period-doubling curve numerically, and proved the local stability changes at $\beta=3/4$ and $1$ by perturbation from the circle \cite{Rob02}. Hu and Sun used the Morse and Maslov-type indices to study the stability of these elliptic triangles \cite{HSindex}. For the Newtonian reduction see Hu, Long, and Sun \cite{HLS}, whose index
theory constructs two analytic branches issuing from \(\beta=3/4\), computes
their circular slopes, and shows that both vanish at collision. Hu, Wu, and Zhu
\cite{HWZ} subsequently proved strict separation and one-dimensionality of the
full antiperiodic kernel at every degeneracy for \(0<e<1\), while the
complex-time argument in \cite[Corollary~B.1]{HWmultiplicity} identifies
precisely two such mass parameters in \([0,9]\) at each \(0<e<1\). Here we
follow both roots globally to the collision endpoint. The homogeneous-potential family studied by Mart\'inez, Sam\`a, and Sim\'o \cite{MSSearly} includes the Newtonian problem treated here; whereas Mart\'inez, Sam\`a, and
Sim\'o traced the stability regions \cite{MSSdiagram} and derived
fixed-stiffness asymptotics under nonvanishing connection hypotheses \cite{MSS}, in the
present problem the stiffness varies with \(e\) to keep the antiperiodic
boundary determinant zero, at a rate that differs between the two reflection
sectors.

\begin{figure}[tbp]
\centering
\includegraphics[width=\linewidth]{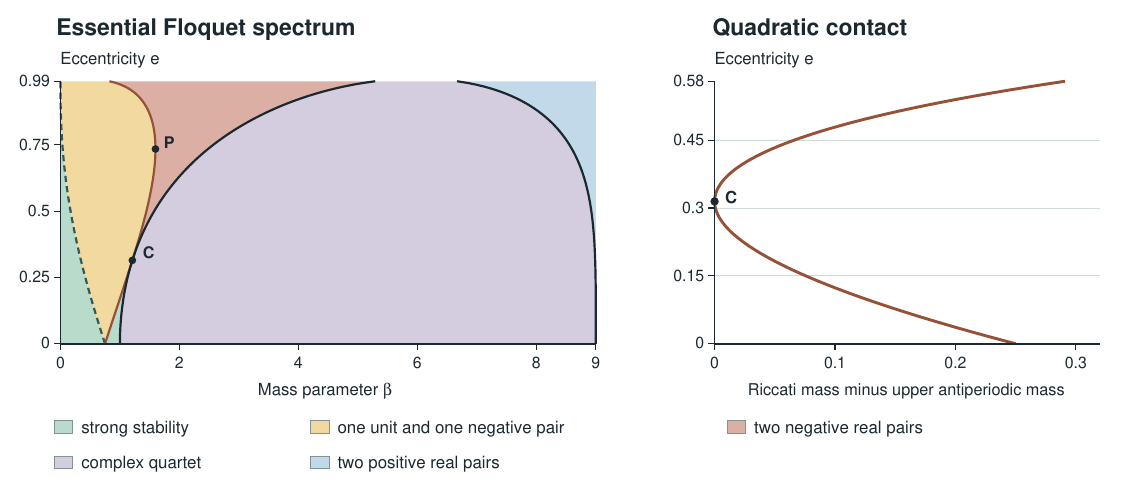}
\caption{Essential spectrum in the \((\beta,e)\)-plane, with \(e\) on the
vertical axis and \(0\le e\le0.99\). The dashed curve is
\(\beta_s\), the solid brown curve is \(\beta_m\), and the black curves are
\(b_-\) and \(b_+\), with the two strong-stability regions sharing the same
shading and \(C,P\) marking the unique quadratic contact and upper-branch
maximum. The right panel plots the horizontal mass gap \(b_- -\beta_m\) against
\(e\) near \(C\). The antiperiodic curves are plotted from truncated Fourier
determinants, while Theorem~\ref{ii-sh-full-diagram} establishes the
classification and boundary blocks.}
\label{ii:full-figure}
\end{figure}

For \(e>0\), denote the lower and upper antiperiodic mass parameters by
\(\beta_s(e)\) and \(\beta_m(e)\), respectively. Our global curve, contact,
collision, and stability-diagram results are summarized in the following three
theorems.
\begin{theorem}[The two physical branches]\label{ii:intro-branches}
There exist real-analytic reflection branches \(\beta_A,\beta_B\) for
\(-1<e<1\), with \(\beta_A(0)=\beta_B(0)=3/4\) and
\(\beta_A'(0)=-\beta_B'(0)=\sqrt{33}/4\), satisfying
\[
 \beta_B(e)=\beta_s(e)<\tfrac34(1-e)<\beta_m(e)=\beta_A(e)
 \qquad(0<e<1).
\]
On \([0,1)\), the lower curve \(\beta_s\) decreases strictly, whereas the upper
curve \(\beta_m\) has a unique stationary point, a nondegenerate maximum at
some \(e_p\in(0,1)\). In particular, for the masses \((40,1,1)\) there
are exactly two antiperiodic eccentricities in \((0,1)\), both lying on the
upper branch and both transverse.
\end{theorem}
\begin{theorem}[The collision endpoint and the hyperbolic boundary]\label{ii:intro-collision}
For \(\rho=(1-e)^{1/4}\) there exist holomorphic germs \(U_m,U_s\) at the
origin, real on real arguments, such that for all sufficiently small
\(\rho>0\),
\[
 \beta_m=\rho U_m(\rho,\rho\log\rho),\qquad
 \beta_s=\rho^{10}U_s(\rho,\rho\log\rho),\qquad
 U_m(0,0)=\frac{12}{2^{1/4}\pi},\quad U_s(0,0)=\frac{8\sqrt2}{\pi^2}.
\]
At a unique \(0<e_c<\sqrt{72/361}<\sqrt{32/153}<e_p\), the upper branch makes quadratic
contact with \(b_-(e)=5+e^2-2\sqrt{2(1-e^2)(2+e^2)}\). On \([0,1)\), the
largest mass parameter \(\beta\) admitting a unit-circle multiplier equals
\(b_-\) for \(e\le e_c\) and \(\beta_m\) for \(e\ge e_c\), and at \(e=e_c\)
the return has a single four-dimensional Jordan block at \(-1\).
\end{theorem}
\begin{theorem}[The global stability diagram]\label{ii:intro-diagram}
With the explicit Riccati arcs
\(b_\pm(e)=5+e^2\pm2\sqrt{2(1-e^2)(2+e^2)}\) of Wang and Zhou \cite{WZ}
(satisfying \(0<\beta_s<\beta_m\le b_-<b_+<9\) on \((0,1)\), with equality only
at \(e=e_c\)), the reduced Floquet discriminant on \((0,9)\times(0,1)\) factors
as \(\Delta=gH^2\), where \(g=(\beta-b_-)(\beta-b_+)/(9-\beta)\) and \(H\) is
real-analytic and nowhere zero. Consequently, the open rectangle
\((0,9)\times(0,1)\) is partitioned into six open regions with simple essential
spectrum:
\begin{itemize}
\item \(0<\beta<\beta_s(e)\) (all \(e\)) and \(\beta_m(e)<\beta<b_-(e)\)
(\(e<e_c\)): two elliptic (nonreal conjugate unit-circle) pairs, giving strong
linear stability under symplectic perturbations;
\item \(\beta_s(e)<\beta<\beta_m(e)\) (all \(e\)): one elliptic pair and one
negative real reciprocal pair;
\item \(\beta_m(e)<\beta<b_-(e)\) (\(e>e_c\)): two negative real reciprocal
pairs;
\item \(b_-(e)<\beta<b_+(e)\) (all \(e\)): a nonreal off-circle quartet;
\item \(b_+(e)<\beta<9\) (all \(e\)): two positive real reciprocal pairs.
\end{itemize}
Moreover, every complex eigenspace on \((0,9)\times(0,1)\) has dimension one:
along \(\beta_s\) (all \(e\)) and \(\beta_m\) (\(e<e_c\)) the spectrum consists
of \(J_2(-1)\) and one elliptic pair; along \(\beta_m\) (\(e>e_c\)) it consists
of \(J_2(-1)\) and one negative real pair; along \(b_-\) it consists of a
conjugate unit-circle pair \(J_2(z),J_2(\bar z)\) for \(e<e_c\) and a negative
real reciprocal pair \(J_2(z),J_2(z^{-1})\) for \(e>e_c\); along \(b_+\) it
consists of a positive real reciprocal pair \(J_2(z),J_2(z^{-1})\); and at
\(\beta_m(e_c)=b_-(e_c)\) it consists of a single block \(J_4(-1)\)
(Theorem~\ref{ii-sh-full-diagram} and Figure~\ref{ii:full-figure}).
\end{theorem}
The three regular boundary edges \(e=0\), \(\beta=0\), and \(\beta=9\) are
likewise classified completely in Section~\ref{ii-sh-diagram} (including the
semisimple double \(-1\) at \((3/4,0)\), the two size-two unit blocks at
\((1,0)\), the Kepler Jordan form \(J_2(1)\oplus J_1(1)^2\) along \(\beta=0\),
and two semisimple double positive reciprocal pairs along \(\beta=9\)).
Moreover, the collision expansion method of Theorem~\ref{ii:intro-collision}
applies beyond three bodies: for every centered regular \(n\)-gon with even
\(n\ge4\), Corollary~\ref{ii-end-polygon} produces an explicit convergent
\((1-e)^{1/4}\) power--log antiperiodic curve in the alternating mode with
positive central mass near collision.

The proofs rest on four analytic ingredients. First, on the upper antiperiodic
branch, eliminating the odd Fourier modes produces a nonnegative quadratic form
whose nullvector variation supplies a global curvature test at every stationary
point (Lemma~\ref{ii:stationary-comparison}); this establishes unimodality
analytically, which together with the endpoint values and the explicit point
bound \(\beta_m(3/5)>243/196\) yields the two transverse crossings for
\((40,1,1)\) (compare also the trigonometric trial lower bound of Wang and Zhou
for the upper branch \cite{WZ}). Second,
along the Riccati arc \(b_-\), monotonicity of the trace on the invariant
Lagrangian graph yields a unique simple crossing of \(-2\), the squared
determinant and nonzero mass derivative make the contact with \(\beta_m\)
quadratic, and the two-block trial bound places \(e_c<e_p\), determining the
branch switch of the hyperbolic boundary. Third, complex-time monodromy
factorization across the two interior singularities \(0\) and \(r\) gives
geometric multiplicity one throughout the open Lagrange domain
(Appendix~\ref{ii-app:simple}, specializing \cite[Theorem~B]{HWmultiplicity}),
ruling out quadratic annihilators and fixing the exact Jordan block sizes on
every interior boundary. Fourth, for the collision asymptotics we use annular
monodromy and exterior powers to isolate the holomorphic boundary pairings,
then apply analytic division and the implicit function theorem in two
independent complex variables before specializing the second variable to
\(\rho\log\rho\).

The remainder of the paper is organized as follows. Section~2 derives the
half-period determinant factorization and identifies the two reflection
sectors for the Lagrange problem. Section~\ref{ii-sh-shape} proves
Theorem~\ref{ii:intro-branches} by establishing global monotonicity of the
lower curve and the uniqueness of the peak on the upper curve,
while Section~\ref{ii-sh-envelope} analyzes the quadratic contact with \(b_-\)
and proves the complete spectral and Jordan classification across every open
region, interior spectral boundary, and regular edge of the parameter rectangle
(Theorem~\ref{ii-sh-full-diagram} and Figure~\ref{ii:full-figure}).
Section~\ref{ii-end-section} constructs the convergent collision expansions of
Theorem~\ref{ii:intro-collision} and applies the same method to the
alternating mode of a regular polygon, and Appendix~\ref{ii-app:simple} gives a
self-contained complex-time proof of geometric simplicity on the open domain.

\section{Half-period determinants and reflection sectors}
Let \(T\) be the half-period map in a fixed symplectic position--velocity
fiber of dimension \(2d\), and suppose that time reversal about \(0\) and
\(\pi\) acts by the same antisymplectic involution \(K\) on that fiber. Its
eigenspaces \(V_\pm\) are complementary Lagrangian planes, reflection of the
second half of the orbit gives \(P=KT^{-1}KT\) for the full return, and we
write \(T_{ab}\) for the block of \(T\) mapping \(V_b\) to \(V_a\).

At the antiperiodic endpoint the characteristic polynomial of the
reversible return reduces to the product of two diagonal minors
\cite[Proposition~3.5]{HWmultiplicity}. We recall the identity with the
forward choice of half-period.

\begin{proposition}[Half-period factorization]\label{ii:reversible}
The antiperiodic nullity is \(\dim\ker T_{++}+\dim\ker T_{--}\), and
\[
 \det(I+P)=2^{2d}\det T_{++}\det T_{--}.
\]
\end{proposition}
\begin{proof}
In bases adapted to the two eigenspaces, we have
\(TK(P+I)=KT+TK=2\operatorname{diag}(T_{++},-T_{--})\). Taking determinants
with \(\det T=1\) and \(\det K=(-1)^d\) gives the determinant formula, while
the kernel identity follows from the same operator relation before taking
determinants.
\end{proof}
Thus a simple zero of one sector determinant \(f(s,e)\) with \(f_s\ne0\)
yields a local analytic mass graph whose derivative is \(-f_e/f_s\), with
second derivative \(-f_{ee}/f_s\) at any stationary point; when both factors
vanish simultaneously, this calculation applies to each factor separately and
the remaining multipliers must still be accounted for.

For \eqref{ii:model}, set \(Q=\operatorname{diag}(1,-1)\) and
\(K=\operatorname{diag}(Q,-Q)\), so that the two mixed problems on
\([0,\pi]\) for \(u=(x,y)\) are
\[
 \begin{array}{c|cc}
 &0&\pi\\ \hline
 A&y=x'=0&x=y'=0\\
 B&x=y'=0&y=x'=0
 \end{array}
\]
in \(L^2((0,\pi);\mathbb R^2)\), where the \(H^2\) operator domain carries
both displayed conditions and the \(H^1\) form domain carries the value
conditions.

At \(\kappa=0\), integrating by parts (using \(xy=0\) at both endpoints) yields
the quadratic form
\[
 \int_0^\pi \bigl[(y'+2x)^2+(x')^2+(3/q-4)x^2\bigr]\,dt .
\]
Because \(w=q\sin t\) is strictly positive on \((0,\pi)\) and solves
\(-w''+(3/q-4)w=0\) with Dirichlet boundary conditions, the scalar Dirichlet
form is nonnegative with kernel spanned by \(w\). Equality in the full
Dirichlet form therefore forces \(x=cw\) and \(y'=-2cw\); since
\(\int_0^\pi w\,dt=2\) and \(y\) vanishes at both endpoints, this implies
\(u=0\), and compact resolvent gives a positive Dirichlet spectral gap.

Consequently, each mixed form domain splits form-orthogonally (by Green's
identity) into the Dirichlet space and the two-dimensional space of homogeneous
solutions extending the prescribed free endpoint values, where the form is
represented by the symmetric matrix mapping those free values to the
corresponding components of \((-p(0),p(\pi))\) with \(p=u'+Ju\). To determine
its inertia, take the homogeneous solution \(x=w\), \(y'=-2w\): in sector \(A\)
with \(y(0)=0\), the free values \((x(0),y(\pi))=(0,-4)\) have outward momenta
\((-(1+e),0)\), while in sector \(B\) with \(y(\pi)=0\), the free values
\((y(0),x(\pi))=(4,0)\) have outward momenta \((0,-(1-e))\). By symmetry of the
boundary matrix, the respective determinants are \(-(1+e)^2/16\) and
\(-(1-e)^2/16\), so each sector has Morse index one and zero nullity for every
\(|e|<1\).

\begin{figure}[htbp]
\centering\includegraphics[width=.91\linewidth]{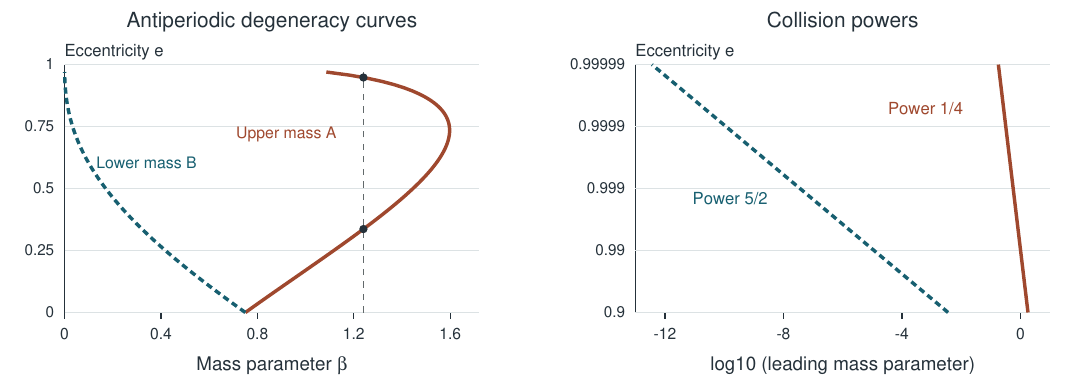}
\caption{Antiperiodic masses of the Lagrange triangle. The left panel shows
the two branches with horizontal mass parameter \(\beta\) and vertical
eccentricity \(e\), up to \(e=0.97\). The vertical line
\(\beta=243/196\) gives the masses \((40,1,1)\). Its intersections with
\(\beta_m\) are marked. The right panel shows the leading terms of
Theorem~\ref{ii:intro-collision}. The ordinate is labeled by \(e\) on
the scale \(-\log_{10}(1-e)\), increasing towards collision.}
\label{ii:branches-figure}
\end{figure}

Put \(K_e=\mathcal A_{3/2,e}\) and \(\Delta=3-2\kappa\). Extending either
sector antiperiodically and rotating via \(u=e^{-Jt}v\), the antiperiodic
Poincar\'e inequality gives \(K_e\ge[3(1-|e|)/(4(1+|e|))]I>0\). Since
\[
 \partial_\kappa(\mathcal A_{\kappa,e}/\Delta)=2K_e/\Delta^2
 \qquad(0\le\kappa<3/2),
\]
the min--max principle combined with the endpoint index counts yields a
unique, simple mass crossing per sector. On its kernel, the derivative of the
original quadratic form is \(2\langle K_e u,u\rangle/\Delta>0\), so the
implicit function theorem on the fixed form domain provides real-analytic
dependence on signed \(e\). At \(e=0\), the zeroth Fourier block in
\eqref{ii-sh-block} decouples and has determinant \(1/5-4\beta/15\); since the
first neighboring block enters the Schur complement only at order \(e^2\), its
first variation in \(e\) is \(\pm2\sqrt{9-\beta}/15\), which gives the circular
slopes in Theorem~\ref{ii:intro-branches}. The chord estimates in
Lemmas~\ref{ii-sh-chord} and~\ref{ii:upper-chord} then distinguish the two
sectors globally. Figure~\ref{ii:branches-figure} shows their two mass
branches and the intersections with the fixed-mass line for $(40,1,1)$.

\section{Global shape of the two curves}\label{ii-sh-shape}

\begin{theorem}[Global monotonicity and unique peak]\label{ii-sh-global}
The lower antiperiodic curve satisfies $\beta_s'(e)<0$ for $0\leq e<1$, whereas
the upper curve $\beta_m$ has a unique stationary point $e_p\in(0,1)$,
increasing strictly on $[0,e_p)$, decreasing strictly on $(e_p,1)$, and
satisfying $\beta_m''(e_p)<0$ (where the derivatives at $e=0$ are right
derivatives).
\end{theorem}

\subsection{An affine Fourier pencil and a lower-branch bound}

Write $k=\kappa$, $d=3-2k$, $p=k(3-k)=\beta/4$, $a=e/2$,
$E=\operatorname{diag}(1,-1)$, and $D=\operatorname{diag}(3-k,k)$, with
forward true anomaly $t$ and $q=1+e\cos t$. On $(0,\pi)$, the upper sector $A$
carries $y(0)=x'(0)=x(\pi)=y'(\pi)=0$, while the lower sector $B$ carries
$x(0)=y'(0)=x'(\pi)=y(\pi)=0$. For $0<e<1$, the two antiperiodic mass
parameters are distinct and each carries a one-dimensional kernel
\cite[Corollary~B.1]{HWmultiplicity}; the mixed boundary conditions identify
the sector at each root, and the mass crossing is positive as shown above
(see also \cite{HLS,HWZ}).

For reference, we recall the Fourier reduction used below \cite[Sec.~4]{HWZ}.
Expand the components in $(\cos\lambda_n t,\sin\lambda_n t)$ in sector~$A$,
and in $(\sin\lambda_n t,-\cos\lambda_n t)$ in sector~$B$, with
$\lambda_n=n+\tfrac12$ and the same normalization in both components, so that
the $n$th block of $T=-\partial_t^2-2J\partial_t$ is
$T_n=\left(\begin{smallmatrix}\lambda_n^2&2\lambda_n\\
2\lambda_n&\lambda_n^2\end{smallmatrix}\right)$. Setting $C=D^{1/2}$, the
bounded pencil $\mathcal L=q+CT^{-1}C$ has the same kernel as the physical
problem under $U=T^{-1}CY$ and $Y=-q^{-1}CU$. Its off-diagonal blocks are $aI$,
its zeroth diagonal block is $M_0+\sigma aE$ with $\sigma=1$ in $A$ and
$\sigma=-1$ in $B$, and the remaining diagonal blocks are
\begin{equation}\label{ii-sh-block}
 M_n=I+
 \begin{pmatrix}c_n(3-k)&b_n\sqrt p\\b_n\sqrt p&c_n k\end{pmatrix},
 \qquad
 c_n=(\lambda_n^2-4)^{-1},\quad b_n=-2c_n/\lambda_n .
\end{equation}
Thus $\mathcal L=M+eV_\sigma$ with $M=\operatorname{diag}_n M_n$, and
differentiating $CT^{-1}C$ on nullvectors yields
$\langle Y,\mathcal L_kY\rangle=-\langle U,q^{-1}D_kU\rangle<0$.

Let $G$ denote the first diagonal block of the inverse of the tail starting at
$n=2$; eliminating this tail and the block $n=1$ yields
\begin{equation}\label{ii-sh-schur}
 R=(M_1-a^2 G)^{-1},\qquad S_\sigma=M_0+\sigma aE-a^2R .
\end{equation}
Since $M_n\geq I$ for $n\geq2$, the tail starting at $n=2$ is bounded below
by $1-e$, whereas the tail starting at $n=1$ is a rank-one negative
perturbation of an operator $\geq(1-e)I$ and has a negative direction because
$\det M_1=-5/7-16p/63<0$. Conjugating by the positive operator's inverse
square root shows that the $n\ge 1$ tail is invertible of index one, so
$M_1-a^2G$ has eigenvalues of opposite signs. At $k=0$, each $S_\sigma$ is
positive: at $e=0$ it equals $\operatorname{diag}(1/5,1)$, and there is no
antiperiodic kernel at $k=0$ for finite $e$ \cite{HLS}. Because Schur
elimination preserves the kernel crossing form (the terms differentiating the
reconstructed vector vanish on the kernel), $S_\sigma$ is positive before its
sole root and has one negative eigenvalue afterward; at the root it has rank
one and positive trace, so $\partial_k\det S_\sigma<0$.

We will also use a finite-dimensional trial for the upper curve, whose root in
\(k\) we denote by \(k_A(e)\). By the preceding inertia decomposition, the full
upper pencil has negative index one for \(0<k\le k_A(e)\), including at the
root, while its restriction to modes \(n=0,1\) has diagonal blocks
\(M_0+aE,M_1\) and off-diagonal block \(aI\). Because \(M_1\) has one negative
eigenvalue, taking its Schur complement yields, for every \(0<k\le3/2\),
\begin{equation}\label{ii-sh-two-block-trial}
 v^T\bigl(M_0+aE-a^2M_1^{-1}\bigr)v<0
 \quad\Longrightarrow\quad \beta_m(e)<4k(3-k),
\end{equation}
since any such \(v\) produces a two-dimensional negative subspace in the
restricted pencil and hence forces \(k>k_A(e)\), while \(4k(3-k)\) increases
strictly on \(0<k\le3/2\).

\begin{lemma}[Lower chord estimate]\label{ii-sh-chord}
For $0<e<1$, $\beta_B(e)<3(1-e)/4$.
\end{lemma}
\begin{proof}
We estimate $S_-$ along the line $p=3(1-e)/16$, where $d^2=(33+3e)/4$. Put
$P_x=\operatorname{diag}(1,0)$; we claim that the full positive tail satisfies
\[
 G\ \geq\frac{10}{7(4-e)}P_x .
\]
Indeed, insert scalar multiples of $(u_n,0)_{n\ge2}$ into the variational
formula for its inverse, where $u_2=1$ and
$u_n=(-1)^{n-2}(21/4)/((n+1/2)^2-1)$ for $n\ge2$. At $(k,e)=(0,1)$, direct
substitution shows that rows $n\geq3$ of the scalar $x$-tail annihilate $u$
and the first row equals $21/10$, while $-\sum_{n\geq2}u_nu_{n+1}=7/10$ by
telescoping reciprocals of four consecutive half-integers. Hence the trial
energy at $(k,e)$ is $21/10+7(1-e)/10-k\sum_{n\ge2}c_nu_n^2\le7(4-e)/10$, and
maximizing over the scalar multiple yields the bound.

Setting $s=5e^2/[14(4-e)]$ and $N=M_1-sP_x$, we have $M_1-a^2G\le N$, and both
matrices have one negative eigenvalue since $\det N=\det M_1-s(1+2d)/7<0$. By
ordered eigenvalue monotonicity, all matrices on the joining segment share the
same inertia, so integrating the inverse derivative yields
$S_-\le U:=M_0-aE-a^2N^{-1}$. With
\[
 F_d=e\left(\frac{128}{315}-\frac{8e}{63}
             -\frac{646e^2}{2205}+\frac{25e^3}{294}\right),
 \qquad d_0=(20+e)/7 ,
\]
expanding the $2\times 2$ determinant and substituting $d^2=(33+3e)/4$ gives
\[
\begin{split}
 (4-e)\det N\,\det U
 ={}&F_d(d-d_0)\\
 &+e(1-e)\left(\frac{1216}{2205}+\frac{208e}{2205}
           -\frac{5906e^2}{15435}+\frac{391e^3}{4116}\right).
\end{split}
\]
Here $d^2-d_0^2=(1-e)(17+4e)/196>0$, while differentiating the cubic $F_d/e$
with $e^2\le e$ shows that $F_d/e$ decreases to $1/14$; since the bracketed
quartic exceeds $2606/15435$, we obtain $\det U<0$. Thus $S_-$ already has a
negative direction on the line, forcing $\beta_B(e)<3(1-e)/4$.
\end{proof}

\begin{lemma}[Upper chord estimate and global sector identification]\label{ii:upper-chord}
For \(0<e<1\), \(\tfrac34(1-e)<\beta_A(e)\). Consequently,
\(\beta_B(e)=\beta_s(e)<\tfrac34(1-e)<\beta_m(e)=\beta_A(e)\) for all
\(0<e<1\).
\end{lemma}
\begin{proof}
Evaluate \eqref{ii-sh-block}--\eqref{ii-sh-schur} at \(p=3(1-e)/16\),
\(d=\sqrt{(33+3e)/4}\), and \(a=e/2\), and set \(r=(1-\sqrt{1-e^2})/2\),
\(\tau=1-r\), so that \(a^2=r\tau\). By the Schur equation \(g=(1-a^2g)^{-1}\)
and continuation from \(a=0\), the free positive half-line matrix with
diagonal one and off-diagonals \(a\) has first inverse entry \(1/\tau\); since
\(M_j\ge I\) for \(j\ge2\), inverse monotonicity gives \(G\le I/\tau\) and
hence \(N:=M_1-rI\le M_1-a^2G\). Because the eigenvalues of \(M_1\) are
\(1/7\pm\sqrt{36/49+112p/441}\), both \(N\) and \(M_1-a^2G\) are invertible of
index one, and their positive and negative eigenvalues stay separated along the
ordered segment joining them. Integrating the derivative of matrix inversion
along this segment therefore gives \(S_+\ge L:=M_0+aE-a^2N^{-1}\). Since the
positive eigenvalue of \(N\) is at least \(\tau\), we have
\(\operatorname{tr} L\ge6/5-a^2/\tau>0\), while expanding \(F=(\det N)\det L\)
yields
\begin{align*}
 F={}&(e/5+4ad/15-a^2)\{\tau(\tau-12/7)-16p/63\}\\
 &-a^2(6\tau/5-12/35+416p/315-4ad/7)+a^4 .
\end{align*}
The coefficient of \(d\) is \(4a\{\tau(3-8\tau)/105-16p/945\}<0\). Since
\(d^2-((20+e)/7)^2=(1-e)(3/28-(1-e)/49)>0\), replacing \(d\) by \((20+e)/7\)
yields an upper bound on \(F\) which, in terms of \(z=\sqrt{(1-e)/(1+e)}\),
factors as
\[
 -\frac{z(1-z)}{8820(1+z^2)^4}
 (4851+15537z+18006z^2+11014z^3+6163z^4+4681z^5+2212z^6)<0.
\]
Since \(\det N<0\), it follows that \(L\) is positive definite, and the sign
change of \(S_+\) at its unique mass root implies \(\beta_A(e)>\tfrac34(1-e)\).
Combining this with Lemma~\ref{ii-sh-chord} yields
\(\beta_B(e)<\tfrac34(1-e)<\beta_A(e)\), identifying \(\beta_B=\beta_s\) and
\(\beta_A=\beta_m\) on \((0,1)\).
\end{proof}
\subsection{Strict decrease of the lower branch}

For $0<e\le1/2$, evaluate $M$ at the lower root, where its only negative
direction lies in block~$1$ because $\det M_0=(3-4\beta)/15>0$ and
$\operatorname{tr}M_0=6/5$. Write $\mathcal L=\mathcal L_B$. If
$\langle\mathcal L u,M^{-1}\mathcal L u\rangle<0$ for some vector $u$, then
$M^{-1}\mathcal L u$ is a negative vector $M$-orthogonal to $\ker\mathcal L$;
conjugating by $|M|^{1/2}$ in its positive and one-dimensional negative summands
and applying Cauchy--Schwarz shows that the $M$-orthogonal complement of such a
vector is positive definite. Hence $\langle Y,MY\rangle>0$ on any nonzero
kernel vector $Y$, so that $\langle Y,V_-Y\rangle=-\langle Y,MY\rangle/e<0$;
differentiating $\mathcal L Y=0$ and using the negative mass crossing then
gives $k_B'<0$. Taking $u_1=(1,0)^T$ with all other blocks zero gives the test
value
\[
 \frac{1-2d}{7}+\frac{e^2}{4}
 \left(\frac{9+2d}{3-4\beta}+\frac{1+4k/9}{7/3+16p/225}\right).
\]
By Lemma~\ref{ii-sh-chord}, the first fraction in parentheses is bounded by
$5/e$, while $d>\sqrt{33}/2>17/6$ and $k<3/8$, so the displayed expression is
less than $-2/3+5e/4+e^2/8\le-1/96$. At $e=0$, the circular slope is
$-\sqrt{33}/4$ \cite{HLS}.

For $1/2<e<1$, rotate $U$ to $u=-JU$, which satisfies
$u''+2Ju'=\operatorname{diag}(k,3-k)u/q$ with lower-sector boundary conditions
$u_2(0)=u_1'(0)=u_1(\pi)=u_2'(\pi)=0$. Following Tschauner \cite{Tschauner,BKV},
set $X=e\cos t$, $B_0=(d-3(2X^2-e^2)/d)/4$, and $c^2=d^2-8+2e^2+9e^4/d^2$,
and define the periodic Riccati matrices
\begin{equation}\label{ii-sh-ric}
 qP_c=\begin{pmatrix}
 q'(1+3X/d)/2 & q+B_0+c/4\\
 -q+B_0-c/4&q'(1-3X/d)/2
 \end{pmatrix}.
\end{equation}
Direct substitution verifies $P_c'+P_c^2+2JP_c=\operatorname{diag}(k,3-k)/q$,
$\operatorname{tr}P_c=q'/q$, and $P_c-P_{-c}=-cJ/(2q)$. For $\beta\le3/4$, let
$C$ be the positive square root in this formula. The graph solutions for
$P_{\pm C}$ with initial position $(1,0)^T$ span the left boundary space
because their initial velocities differ, while the right boundary columns are
their first positions at $\pi$ times $(1,(P_{\pm C})_{21}(\pi))^T$, which are
linearly independent. At $e=0$ and $\beta=3/4$, the lower-sector kernel
$(\cos(t/2),-(7-\sqrt{33})\sin(t/2)/4)^T$ has initial velocity belonging to
$P_{-C}$ (with $C=1/2$). Since the sector kernel is simple, the two scalar
first-position factors cannot vanish simultaneously along the branch; by
connectedness and Lemma~\ref{ii-sh-chord}, the $P_{-C}$ graph supplies the
kernel for all $0\le e<1$.

Fix $1/2<e<1$ and $0\le k\le k_B(e)$, and set
\[
 m=d-3e^2/d,\quad g=m+C,\quad N=q(P_{-C})_{12},\quad L=2q-N .
\]
Both $N$ and $L$ are positive: indeed, $d>d_0=(20+e)/7$ and $d>3e$ by
Lemma~\ref{ii-sh-chord}, and the identity
$d_0^2-(3-e)d_0-3e^2=(1-e)(139e-20)/49>0$ gives $m>3-e$. Moreover,
$m^2-C^2=8(1-e^2)$; viewed as a concave quadratic in $X$, $N$ has minimum
$1-e+(m-C)/4>0$, while $L_X=1+3X/d>0$ gives
$\min L=(1-e)(1-2(1+e)/g)>0$ (when $m<2(1+e)$, use
$C^2-[2(1+e)-m]^2=4(1+e)(m-3+e)>0$).

The identity $q(P_{-C})_{11}=L'/2$ prompts the substitution $u_1=Lz$ and
$u_2=(qu_1'-L'u_1/2)/N$, which gives $u_1'-2u_2=-(L^2/N)z'$ and converts the
first physical equation into the Sturm--Liouville problem
\begin{equation}\label{ii-sh-sturm}
 -(Pz')'=kwz,\qquad P=L^2/N,\quad w=L/q,\quad z'(0)=z(\pi)=0 .
\end{equation}
Conversely, the first graph residual $v=u'-P_{-C}u$ has $v_1=0$, and the first
component of the physical residual $v'+(P_{-C}+2J)v$ equals $-(L/q)v_2$, so
\eqref{ii-sh-sturm} implies both graph equations while the vanishing diagonal
entries at $0,\pi$ enforce the mixed boundary conditions. On the fixed $H^1$
form domain $z(\pi)=0$, consider $Q[z]=\int_0^\pi(Pz'^2-kwz^2)\,dt$. It is
positive at $k=0$ and remains nondegenerate for $0\le k<k_B(e)$ by the
reconstruction and sector root count, so continuity of the lowest eigenvalue
makes the zero at $k=k_B(e)>0$ a ground state. Normalizing $z(0)=1$ and setting
$F(t)=\int_0^t wz$, we have $z>0$ on $[0,\pi)$ and $Pz'=-kF$, so $F$ increases
and $z$ decreases strictly on $(0,\pi)$.

We now determine the sign of $Q_e[z]$ at this kernel, differentiating always at
fixed $k$. Put $S=1-e^2$, $b=3/(2d)$, $K=2/g-b$, and $H=SK$, so that
expansion in $q$ yields
\[
 L=bq^2-(2k/d)q-H,\qquad w=bq-2k/d-H/q .
\]
Since $(d/3+3e^2/d)^2-C^2=8\beta/9>0$ at the root, we have $4d/3-g>0$ and
$g_e=6e(4d/3-g)/(dC)>0$, whence $K>0$ and $K_e<0$. Using
$S\cos t/q^2=(\sin t/q)'-e/q$ then yields
\begin{equation}\label{ii-sh-prefix}
 \int_0^t w_e
 = (b+K/q(t))\sin t+(eK-SK_e)\int_0^t q^{-1}>0
 \qquad(0<t\le\pi).
\end{equation}

For the reciprocal coefficient $r=P^{-1}$, the required inequality runs from the
opposite endpoint. Setting $\alpha=4k/3$ and $\gamma=H/b$, so that
$L/b=q^2-\alpha q-\gamma$, we have
\[
 -\gamma_e/\gamma=2e/S+8e/(Cg)<4e/S,
\]
because $C^2=1-\beta+2e^2+9e^4/d^2>9/8$ by the chord bound and $e>1/2$,
giving $Cg>3\ge4S$. Since $\alpha/q+\gamma/q^2<1$, expanding
$r=(2q-L)/L^2$ gives
\[
 r=b^{-2}\sum_{i,j\ge0}\binom{i+j}{j}\alpha^i\gamma^j
 q^{-m_{ij}}(A_{ij}-bq),
\]
where $m_{ij}=3+i+2j$ and $A_{ij}=2(i+j+1)$. Grouping terms with equal $i+j$
and using the strict geometric bound on compact $t$-intervals shows that this
series and its $e$-derivative converge locally uniformly in $(t,e)$, with every
summand positive. Writing $m_1=m_{ij}$ and $A_1=A_{ij}$, the logarithmic
$e$-derivative of a summand is
\[
 j\gamma_e/\gamma-m_1\cos t/q-b\cos t/(A_1-bq),
\]
which increases strictly with $t$ because $A_1>b$ and both rational functions of
$\cos t$ have strictly negative derivative after including their signs.

For $I_m=\int_0^\pi q^{-m}dt$, substituting $\cos t=(\cos s-e)/(1-e\cos s)$
yields $I_m=S^{1/2-m}\int_0^\pi(1-e\cos s)^{m-1}ds$. Since the latter integral
has nonnegative coefficients in $e^2$, we have $I_m'/I_m\ge(2m-1)e/S$;
moreover, $1/q$ and $-\cos t/q$ are similarly ordered functions of $\cos t$,
so their nonnegative covariance and nonnegative expectation of $-\cos t/q$ under
the weight $q^{1-m}$ give $I_m'/I_m\ge I_{m-1}'/I_{m-1}$ for $m\ge3$. Applying
these two moment inequalities to $A_1 I_{m_1}-bI_{m_1-1}$, the logarithmic
derivative of the integral of each summand is at least
\[
 j\gamma_e/\gamma+I_{m_1}'/I_{m_1}\ \ge\ (5+2i)e/S>0.
\]
Because the derivative of each summand can change sign only from negative to
positive while its full integral is positive, every tail integral of that
derivative is positive, and summation yields $\int_t^\pi r_e>0$ for
$0\le t<\pi$.

Integrating the positive tail integral against $(F^2)'$ in the first term and
\eqref{ii-sh-prefix} against $-(z^2)'$ in the second then gives
\[
 Q_e[z]=-k^2\int_0^\pi r_eF^2\,dt-k\int_0^\pi w_e z^2\,dt<0 .
\]
At fixed $e$, positivity below the root implies $Q_k[z]\le0$ via the left
difference quotient on this same $z$. Since the solution normalized by
$z(0)=1$ depends analytically on $e$ along the branch and its derivative belongs to the
fixed form domain (where it pairs to zero with the kernel), differentiating
$Q[z]=0$ along the branch yields $Q_k[z]k_B'+Q_e[z]=0$, whence $Q_k[z]<0$ and
$k_B'<0$. Multiplying by $4d>0$ proves the lower-branch assertion of
Theorem~\ref{ii-sh-global}.

\subsection{Curvature of the upper branch}

Write \(a=e/2\) and \(\widehat k(a)=k_A(2a)\). In the upper pencil
\(\mathcal L\) of \eqref{ii-sh-block}, each off-diagonal block joins an even
mode to an odd mode; grouping the even modes first with \(\Pi y=y_0\)
and \((\mathcal Bv)_{2j}=v_{2j-1}+v_{2j+1}\) (where \(v_{-1}=0\)), the pencil takes the block form
\[
 \begin{pmatrix}
 \mathcal D_{\rm ev}+a\Pi^*E\Pi&a\mathcal B\\
 a\mathcal B^*&\mathcal D_{\rm od}
 \end{pmatrix},
 \qquad \mathcal D_{\rm ev}=\bigoplus_{j\ge0}M_{2j},\quad
 \mathcal D_{\rm od}=\bigoplus_{j\ge0}M_{2j+1}.
\]
Since the odd diagonal is invertible with exactly one negative eigenvalue,
completing the square on all odd modes yields the even-mode pencil
\begin{equation}\label{ii:even-pencil}
 K(k,a)=\mathcal D_{\rm ev}+a\Pi^*E\Pi-a^2\mathcal H,\qquad
 \mathcal H=\mathcal B\mathcal D_{\rm od}^{-1}\mathcal B^* .
\end{equation}
Along \(k=\widehat k(a)\), the full pencil has one negative direction and
a simple kernel, so \(K\) is positive semidefinite with a simple kernel;
a nullvector \(y\) lifts to the full pencil by adjoining
\(-a\mathcal D_{\rm od}^{-1}\mathcal B^*y\), whence \(\langle y,K_k y\rangle<0\)
by the mass crossing in \eqref{ii-sh-block}. Moreover, \(y_0\ne0\) whenever \(y\ne0\),
since otherwise the zeroth row and successive rows of the full equation would
force every component to vanish. Because only the zeroth block carries a term
linear in \(a\), we test the second variation by perturbing the component on
which this linear term is negative.

\begin{lemma}[A curvature test at a stationary mass]\label{ii:stationary-comparison}
Consider a pencil of the form \eqref{ii:even-pencil} on a real Hilbert
space, with \(\Pi\) the projection onto a two-dimensional coordinate block,
bounded self-adjoint \(\mathcal D_{\rm ev},\mathcal H\) depending smoothly on \(k\),
\(\Pi\mathcal D_{\rm ev}=M_0\Pi\), and \(E=\operatorname{diag}(1,-1)\).
Let \(k(a)\) be a twice differentiable curve along which \(K\ge0\).
At a stationary point \(a>0\) suppose that \(Ky=0\), \(u=\Pi y\ne0\),
and \(\langle y,K_k y\rangle<0\). With \(t=(0,1)^T\), put
\[
 c=\langle\Pi^*t,K\Pi^*t\rangle,\qquad
 g=(2M_0/a+E)t .
\]
If \(c>0\) and \(2a(g_2^2-g_1^2)>c\), then \(k''(a)<0\).
\end{lemma}
\begin{proof}
At a stationary point, nonnegativity forces the derivative of
\(\langle y,K(k(a),a)y\rangle\) to vanish, so
\(2a\langle y,\mathcal H y\rangle=u^TEu\), while \(Ky=0\) gives
\(K_a y=-(2\mathcal D_{\rm ev}/a+\Pi^*E\Pi)y\).
For real \(\delta\), evaluating \(K(k(a+\delta),a+\delta)\) on
\(y+\delta\Pi^*t(g^Tu)/c\) yields a nonnegative quadratic form vanishing to
first order at \(\delta=0\), whose second derivative implies
\[
 0\le k''(a)\langle y,K_k y\rangle
       -\frac1a\,u^T\!\left(E+\frac{2a}{c}gg^T\right)u .
\]
Since the matrix in parentheses has positive first diagonal entry and determinant
\(2a(g_2^2-g_1^2)/c-1>0\), the strict mass crossing yields \(k''(a)<0\).
\end{proof}

We now verify the hypotheses of Lemma~\ref{ii:stationary-comparison} uniformly,
beginning with a mass upper bound. The zeroth diagonal block of \(K\) is
\(U=M_0+aE-a^2 M_1^{-1}\); at \(k=3/20\), direct substitution in
\eqref{ii-sh-block} with \(v=a-2/5\) gives
\[
 \det U=-\frac1{375}+\frac{41}{450}v
 -v^2\left(\frac{121}{45}+\frac{275}{72}v+\frac{175}{144}v^2\right).
\]
Because \(v\ge-2/5\), the bracket is at least \(209/180>1\), so completing the
square yields \(\det U\le-1/375+(41/900)^2<0\), and the two-mode trial
\eqref{ii-sh-two-block-trial} gives \(\widehat k(a)<3/20\) for all \(0<a<1/2\).

Write \(b=(M_0)_{12}=16\sqrt p/15\), \(z=(M_0)_{22}=1-4k/15\), and
\[
 r=(M_1^{-1})_{22}=\frac{9(5-4k)}{45+16p}.
\]
On \(0\le k\le3/20\), \(r\) is positive and decreasing with \(r\le1\) and
\(z\ge24/25\), so \(c=z-a-a^2r\ge21/100\); with \(g\) as in the lemma,
\[
 a\{c-2a(g_2^2-g_1^2)\}
   =8(b^2-z^2)+9az-3a^2-a^3r .
\]
The right-hand side increases in both \(a\) and \(k\) on the closed rectangle
\(0\le a\le1/2\), \(0\le k\le3/20\), since its \(a\)-derivative is at least
\(9(24/25)-3-3/4>0\) and its \(k\)-derivative is at least
\(8(256/225)(27/10)-6/5>0\). At the upper corner we have
\(b^2=304/625\), \(z=24/25\), and \(r=55/72\), giving \(-2551/360000<0\), so
Lemma~\ref{ii:stationary-comparison} yields \(\widehat k''<0\) and hence
\(\beta_m''(e)=(3-2k)\widehat k''(a)<0\) at every stationary point.

Because the upper curve has positive right derivative at \(e=0\) and tends to
zero at collision by \cite{HLS} (and independently by Theorem~\ref{ii-end-theorem}),
it attains a maximum in \((0,1)\). Every stationary point is a strict downward
crossing of \(\beta_m'\), leaving a unique stationary point and proving the
shape assertion of Theorem~\ref{ii-sh-global}.

For later geometric comparisons we record a rational lower bound at
\(k=16/125\), \(e=3/5\): the free-tail bound in the proof of
Lemma~\ref{ii:upper-chord} gives \(a^2G\le I/10\), and since the two matrices
being inverted have index one throughout the ordered segment between them,
\[
 S_+\ \ge\ L=M_0+3E/10-(9/100)(M_1-I/10)^{-1}.
\]
Substitution in \eqref{ii-sh-block} gives \(L_{11}>3/5\) and \(\det L>9/1000\), whence
\begin{equation}\label{ii:peak-lower}
 \beta_m(3/5)>4(16/125)(3-16/125)=22976/15625 .
\end{equation}
(Numerical evaluation of the truncated Schur determinant places the unique peak
near \(e_p\approx0.734543\) with \(\beta_m(e_p)\approx1.597496\).)

\begin{corollary}[Transverse intersections for masses $(40,1,1)$]\label{ii-sh-masses}
The masses $(40,1,1)$ admit exactly two antiperiodic degeneracy eccentricities
in $(0,1)$, both lying on the upper curve and both transverse.
\end{corollary}
\begin{proof}
The mass parameter $27(40+40+1)/42^2=243/196$ lies strictly between the
circular value $3/4$ and the lower bound for $\beta_m(3/5)$ in
\eqref{ii:peak-lower}. Since $\beta_m(e)\to0$ at collision,
Theorem~\ref{ii-sh-global} gives one intersection with this mass level
on each side of $e_p$ with nonzero derivatives of opposite signs, while
the lower curve remains strictly below $3/4$.
\end{proof}

\section{Quadratic contact and the hyperbolic boundary}\label{ii-sh-envelope}

Let $M(\beta,e)$ be the full-period essential return, and for each $e\in[0,1)$
let $\beta_k(e)$ be the largest physical mass parameter at which $M$ has a
multiplier on the unit circle, as in Hu--Long--Sun \cite{HLS}; sufficient
conditions and quantitative estimates for the stable and hyperbolic regions in
the $(\beta,e)$-plane were further obtained by Hu and Ou \cite{HOregion}, Hu,
Ou, and Wang \cite{HOWtrace}, and Hu, Ou, and Sun \cite{HOS}. To identify the
boundary itself, define the explicit collision curves
\[
 b_\pm(e)=5+e^2\pm2\sqrt{2(1-e^2)(2+e^2)} .
\]
Wang and Zhou \cite{WZ} derived these arcs via Riccati decomposition, showing
that the nonreal off-circle quartet lies between them; here we determine where
the lower explicit arc $b_-$ meets the antiperiodic curve $\beta_m$ and which
curve supplies $\beta_k$ on each side of the contact.

\begin{proposition}[Strict trace monotonicity on the lower Riccati arc]
\label{ii-sh-arc-trace}
In the rotated coordinates \(u=-JU\) of \eqref{ii-sh-ric}, let \(h(e)\) be the trace of the return on the invariant Lagrangian
graph \((u'+Ju,u)=((P_0+J)u,u)\) at \(\beta=b_-(e)\). Then \(h'(e)<0\) for every
\(0<e<1\), and the function \(T_{\mathrm{arc}}(e)=4-2h(e)\) tends to
\(4-4\cos(\sqrt{2}\pi)\) as \(e\downarrow0\) and to \(+\infty\) as
\(e\uparrow1\). In the full two-factor normalization
\(M=(I+iA)(I+iB)\) of \cite[Lemma~5.1]{HWmultiplicity}, with \(A^2=B^2=0\), one has
\(T_{\mathrm{arc}}=\operatorname{tr}(AB)\).
\end{proposition}

\begin{proof}
Along either explicit arc the Riccati parameter in \eqref{ii-sh-ric} vanishes,
because its square equals \((\beta-b_-)(\beta-b_+)/(9-\beta)\); the Riccati
identity makes the graph invariant, \(P_0-P_0^T=-2J\) makes its canonical lift
Lagrangian, and \(\operatorname{tr}P_0=q'/q\) gives the graph return
determinant one. The symplectic pairing identifies the quotient return with the
inverse transpose of the graph return, so \(\operatorname{tr}M=2h\).

Put \(a=e/(1+\sqrt{1-e^2})\), \(\sigma=a^2\), and \(z=\exp(it)=-ay\), so that
the inner pole is \(y=1\), the outer pole is \(y=1/\sigma\), and
\(e=2\sqrt{\sigma}/(1+\sigma)\). Define
\[
 p_-=\frac{3e^2\sigma}{d^2},\qquad p_+=\frac{d^2\sigma}{3e^2},
 \qquad
 \ell=\frac{2+\sigma+\sqrt{1+4\sigma+\sigma^2}}3 .
\]
On \(b_-\) the arc equation gives \(p_-p_+=\sigma^2\) and
\(p_-+p_+=2(1+\sigma+\sigma^2)/3\); since
\(d^2-3e^2=4(1-e^2)+2\sqrt{2(1-e^2)(2+e^2)}>0\), we have \(p_-<p_+\) and hence
\[
 p_-=1-(1-\sigma)\ell,\qquad
 p_+=1-\frac{1-\sigma}{3\ell}.
\]
In the graph equation \(u'=P_0u\), set \(u=(x_1,ix_2)^T\) and \(\xi=-z\), pass
to the coordinates \(x_1+x_2\) and \(d\xi(x_1-x_2)/(2e)\), and divide both by
\(\sqrt{\xi(-q)}\) (taking the positive square root on \(0<y<1\)). Because
\(X=e\cos t=-e(\xi+\xi^{-1})/2\) and \(Y=zq_z=e(\xi^{-1}-\xi)/2\), the upper
and lower numerators in these sum and difference coordinates are
\(-(d-3(X+Y)^2/d)/4\) and \((d-3(X-Y)^2/d)/4\), respectively; substituting
\(\xi(-q)=\frac e2(1-y)(1-\sigma y)\) produces the traceless system
\begin{equation}\label{ii-sh-arc-connection}
 \frac{dW}{dy}=
 \begin{pmatrix}
 \dfrac1{2y} & \dfrac{1-p_-y^2}{y(1-y)(1-\sigma y)}\\[4pt]
 \dfrac{3(1-p_+y^2)}{4y(1-y)(1-\sigma y)} & -\dfrac1{2y}
 \end{pmatrix}W .
\end{equation}

At \(y=0\) the residue is
\(R_0=\left(\begin{smallmatrix}1/2&1\\3/4&-1/2\end{smallmatrix}\right)\)
with exponents \(1,-1\) and eigenvectors \(v_+=(1,1/2)^T\),
\(v_-=(1,-3/2)^T\), and we write \(W_0=y(v_++O(y))\) for the high-exponent
solution. At \(y=1\), in the coordinate \(s=1-y\), the residue is
\(-\left(\begin{smallmatrix}0&\ell\\1/(4\ell)&0\end{smallmatrix}\right)\) with
exponents \(-1/2,1/2\) and eigenvectors \(w_\pm=(1,\pm1/(2\ell))^T\); writing
\(W_r=s^{1/2}(w_-+O(s))\) for the high-exponent solution there, we set
\(\gamma(\sigma)=\lim_{y\uparrow1}s^{1/2}(W_0)_1\).

In the low-exponent solutions, the coefficients multiplying the high-exponent
solution times the local logarithm are
\begin{equation}\label{ii-sh-arc-resonance}
 n_0=\frac{1+10\sigma+\sigma^2}{8}+\tfrac32 p_-,
 \qquad
 n_r=\frac{3\ell-\ell^{-1}-1+\sigma}{2(1-\sigma)} .
\end{equation}
For the first identity, expand the coefficient matrix as
\(R_0/y+S_0+yS_1+O(y^2)\), where \(S_0,S_1\) have zero diagonal, \(S_0\) has
off-diagonal entries \(1+\sigma\) and \(3(1+\sigma)/4\), and \(S_1\) has
off-diagonal entries \(1+\sigma+\sigma^2-p_-\) and
\(3(1+\sigma+\sigma^2-p_+)/4\); because \(R_0^2=I\), the order-one coefficient
of the low-exponent solution is \(-R_0S_0v_-\) and the order-two resonant
projection is \((3/4,1/2)(S_1-S_0R_0S_0)v_-\). For the second identity, the
constant term at \(s=0\) is
\[
 S_r=\begin{pmatrix}
 -1/2 & -(2-\ell)/(1-\sigma)\\
 -(3/2-1/(4\ell))/(1-\sigma) & 1/2
 \end{pmatrix},
\]
whose order-one resonant projection \((1/2,-\ell)S_rw_+\) yields the second
formula in \eqref{ii-sh-arc-resonance}.

The Wronskian \(\det(W_0,W_r)\) is constant and equals \(-\gamma/\ell\) at
\(y=1\); since \(\det(v_+,v_-)=-2\), the low-exponent coefficient of \(W_r\) at
\(0\) is \(\gamma/(2\ell)\). The two rank-one logarithmic operators \(N_0,N_r\)
send the corresponding low-exponent coefficients to \(n_0W_0\) and \(n_rW_r\),
giving \(\operatorname{tr}(N_0N_r)=n_0n_r\gamma^2/(2\ell)\). Because the
forward physical circle \(|y|=1/a\) encloses \(0\) and \(1\) while excluding
\(1/\sigma\), its monodromy is conjugate to the product of the two positively
oriented local monodromies based at the same point of \((0,1)\). In
\eqref{ii-sh-arc-connection} these factors are \(I+2\pi iN_0\) and
\(-(I+2\pi iN_r)\), and undoing the square-root gauge removes the minus sign
because \(\xi(-q)\) has a single zero inside the physical loop, yielding
\begin{equation}\label{ii-sh-arc-trace-product}
 h=2-2\pi^2 n_0(n_r/\ell)\gamma^2,\qquad
 T_{\mathrm{arc}}=4\pi^2 n_0(n_r/\ell)\gamma^2 ,
\end{equation}
while expanding the two square-zero factors of the full return gives
\(\operatorname{tr}M=4-\operatorname{tr}(AB)\).

To compare \(\gamma\), multiply the second coordinate of \(W\) by
\(L(y)=1+(\ell-1)y\), which transforms the off-diagonal entries to
\[
 \frac1{y(1-y)}+\frac{1+\sigma-\ell}{(1-\sigma y)L(y)}
 \quad\text{and}\quad
 \frac{3}{4y}+\frac1{4(1-y)}+\frac{\ell-1}{4}(3-\ell^{-1})
 +\left(\frac{3(1+\sigma)}4-\frac1{4\ell}\right)
   \frac{L(y)}{1-\sigma y},
\]
while leaving the diagonal entries unchanged except for adding
\((\ell-1)/L(y)\) to the second diagonal entry. From
\(\sigma=(3\ell^2-4\ell+1)/(2\ell)\) we obtain
\(d\sigma/d\ell=(3-\ell^{-2})/2\); since \(1<\ell\) and \(0<\ell'<1\) for
\(0<\sigma<1\), all parameter-dependent terms are nonnegative and nondecreasing
in \(\sigma\) (using \(1+\sigma-\ell>0\) and
\(\partial_\sigma((1-\sigma y)L(y))<0\) in the first off-diagonal entry).

Writing the transformed system as \(V'=(R_0/y+E_\sigma)V\), the matrix
\(E_\sigma\) is bounded, entrywise nonnegative, and nondecreasing in \(\sigma\)
on every \([0,b]\) with \(b<1\). The Volterra equation for \(V/y\) starts at
\(v_+\) with kernel \((t/y)(y/t)^{R_0}E_\sigma(t)\) (\(0<t\le y\le b\)), which
is bounded (as \(R_0\) has exponents \(1,-1\)) and nonnegative (as \(R_0\) has
positive off-diagonal entries), so its convergent iteration yields positive
components that are nondecreasing in \(\sigma\). Because the first component is
unchanged by the gauge, \(\gamma\) is nonnegative and nondecreasing; moreover,
\(\gamma=0\) would force \(W_0\) to be a nonzero multiple of the high-exponent
solution at \(1\), whose components have opposite signs near \(1\), whence
\(\gamma>0\).

Both \(n_0\) (using \(p_-'=\ell-(1-\sigma)\ell'>0\)) and
\[
 n_r/\ell=\frac{3-\ell^{-2}}{2(1-\sigma)}-\frac1{2\ell}
\]
are positive with positive derivatives on \(0<\sigma<1\). Since the physical
trace is analytic on the open arc, taking right difference quotients in
\eqref{ii-sh-arc-trace-product} and using
\(\gamma(\sigma+\epsilon)\ge \gamma(\sigma)>0\) yields
\(\partial_\sigma T_{\mathrm{arc}}>0\) and thus \(h'(e)<0\). Finally, the
circular trace follows from \(P_0^2=-I/2\) at \(e=0\), while as \(\sigma\uparrow1\)
we have \(n_r/\ell\to+\infty\) with the remaining two factors bounded below
by positive constants on every \([\sigma_*,1)\) (\(\sigma_*>0\)).
\end{proof}

\begin{theorem}[Quadratic contact and exchange of boundary branches]\label{ii-sh-contact}
There is a unique $e_c\in(0,1)$ such that $b_-(e_c)=\beta_m(e_c)$, and it
satisfies $e_c<\sqrt{72/361}<\sqrt{32/153}<e_p$. At $(\beta_m(e_c),e_c)$, the return $M$
has a single Jordan block of size four at $-1$, and
\[
 b_-(e)-\beta_m(e)=K(e-e_c)^2+O((e-e_c)^3),\qquad K>0 .
\]
Moreover, $\beta_k=b_-$ on $[0,e_c]$ and $\beta_k=\beta_m$ on $[e_c,1)$, so
that $\beta_k$ is $C^1$ at $e_c$ with right second derivative smaller than its
left second derivative by $2K$.
\end{theorem}

\begin{proof}
On either explicit arc, the determinant-one graph return and its dual quotient give
\[
  \det(\lambda I-M)=(\lambda^2-h\lambda+1)^2,\qquad
  F:=\det(M+I)=(h+2)^2 .
\]
Integrating by parts with periodic $P_0$ and substituting its
Riccati equation shows that the original operator in these coordinates equals
$(\partial_t-P_0)^*(\partial_t-P_0)$ on every unit-twisted form domain. In
particular, the antiperiodic form is nonnegative along $b_-$. At fixed $e>0$,
since the two physical mass crossings are simple and positive and the form is
positive at $\beta=9$, the negative index takes the values $2$, $1$, and $0$
on the three consecutive intervals separated by $\beta_s<\beta_m$
\cite{HLS,HWZ}, forcing $b_-\ge\beta_m$.

By Proposition~\ref{ii-sh-arc-trace}, \(h\) decreases strictly from
\(2\cos(\sqrt2\pi)\in(-2,0)\) to \(-\infty\), yielding a unique \(e_c\in(0,1)\)
with \(h(e_c)=-2\) and \(h'(e_c)<0\); since \(F=(h+2)^2\) on \(b_-\), this is
the sole antiperiodic point on the lower arc, and the index comparison above
identifies its mass with \(\beta_m(e_c)\).

For the upper bound on \(e_c\), set \(G(\sigma)=n_0 n_r/\ell\). At
\(\sigma=1/18\), the relations \(p_-p_+=\sigma^2\) and \(p_+<1\) give
\(p_->\sigma^2\), while the explicit formula for \(\ell\) gives \(\ell>21/20\);
using \(n_0(0)=1/8\) and \((n_r/\ell)(0)=1/2\), we obtain
\[
 \frac{G(1/18)}{G(0)}>
 \frac{517}{324}\,\frac{3158}{2499}>2.
\]
By \eqref{ii-sh-arc-trace-product}, \(T_{\mathrm{arc}}(e)/G(\sigma)=4\pi^2\gamma(\sigma)^2\)
is nondecreasing, so comparison with \(T_{\mathrm{arc}}(0)=4-4\cos(\sqrt2\pi)>4\)
gives \(T_{\mathrm{arc}}(e)>8\) at \(e^2=72/361\) and hence
\(e_c<\sqrt{72/361}<e_0:=\sqrt{32/153}\).

To compare \(e_0\) with the upper-branch maximum, substitute \(k=16/125\),
\(p=k(3-k)\), and \(v=(-9\sqrt p/5,1)^T\) into \eqref{ii-sh-two-block-trial},
which by \eqref{ii-sh-block} yields the convex quadratic
\[
 Q(e)=-\frac{122247229}{732421875}+\frac{74639}{781250}\,e
       +\frac{4525456311}{9937862500}\,e^2
\]
satisfying \(Q(9/20),Q(1/2)<-1/200\). Since \(9/20<e_0<1/2\), the trial
inequality and \eqref{ii:peak-lower} give
\(\beta_m(e_0)<22976/15625<\beta_m(3/5)\), and because \(e_0<3/5\) while
\(\beta_m\) decreases strictly past its unique maximum, it follows that
\(e_0<e_p\).

\smallskip
\noindent\emph{Order of contact and choice of boundary.}
At the contact $e_c$, the full antiperiodic kernel is simple and
$\det(\lambda I-M)=(\lambda+1)^4$, proving that $M$ has a single Jordan block
of size four at $-1$. To show that $F_\beta\ne0$, let $v$ span $\ker(M+I)$ in
canonical coordinates; because $M$ is symplectic and $M+I$ has rank three,
$\operatorname{adj}(M+I)$ is a nonzero multiple of $vv^TJ_4$, where $J_4$ is
the canonical symplectic matrix. Pairing the variation of constants identity
with the solution starting at $v$ yields
$|v^TJ_4 M_\beta v|=|\langle u,q^{-1}(\partial_\beta\operatorname{diag}(k,3-k))u\rangle|>0$
by the positive physical mass crossing, whence Jacobi's determinant formula
gives $F_\beta\ne0$. On $b_-$ we have $F=(h+2)^2$, whereas
$F(\beta_m(e),e)=0$, so Taylor expansion in $\beta$ yields
$K=h'(e_c)^2/F_\beta(\beta_m(e_c),e_c)$; nonnegativity of $b_--\beta_m$
then forces $K>0$. Note also that a quadruple multiplier $-1$ cannot occur
anywhere else: if $c\ne0$, the two complementary graphs in \eqref{ii-sh-ric}
(over $\mathbb C$ if necessary) would each contribute a $-1$ eigenvector, while
on $b_+$ one has $h>2$ by continuity from $(9,0)$ (where
$\operatorname{tr}M/2=2\cosh(\sqrt2\pi)$) and by periodic nondegeneracy
forbidding $h=2$.

At $\beta=9$, periodic rotation removes the gyroscopic term and yields the
scalar quadratic form
\[
 \int_0^{2\pi}\!\left(|v'|^2+
       \left(\frac{3}{2q}-1\right)|v|^2\right)dt
 =\int_0^{2\pi}\!\left(q^2|(v/q)'|^2+\frac{|v|^2}{2q}\right)dt
\]
in each component, since $q''=1-q$ and the boundary terms cancel for every
unit twist. For any $\beta>b_-$, the physical form is a convex combination of
the forms at $b_-$ and $9$ (with coefficient $\sqrt{(9-\beta)/(9-b_-)}<1$
multiplying the nonnegative form at $b_-$) and is therefore strictly positive
for all unit twists. On $b_-$, we have $h(0)=2\cos(\sqrt2\pi)$ and $h$ cannot
reach $2$, so the unique downward crossing of $-2$ at $e_c$ leaves unit-circle
spectrum on $b_-$ precisely for $0\le e\le e_c$. For $e>e_c$ and
$\beta_m<\beta<b_-$, the two real graph returns for $c$ and $-c$ have
determinant one and traces strictly below $-2$ just to the left of $b_-$ by
continuity from $P_0$; neither trace can reach $-2$ or $2$ on this open mass
interval by the antiperiodic root count and periodic nondegeneracy. Hence the
largest mass with unit-circle spectrum is $\beta_m$ for $e\ge e_c$, and the
$C^1$ regularity and jump in one-sided second derivatives at $e_c$ follow from
the quadratic gap.
\end{proof}
\subsection{The stability diagram of the elliptic Lagrange family}
\label{ii-sh-diagram}

Write the characteristic polynomial of the essential return as
$z^4-a z^3+b z^2-a z+1$, so that the two roots of $h^2-ah+b-2=0$ are the
reduced traces $h=z+z^{-1}$ with discriminant $\Delta=a^2-4(b-2)$. Below, an
elliptic pair means two nonreal conjugate unit multipliers, while a positive or
negative pair consists of distinct real reciprocal multipliers of the indicated
sign. The elliptic and mixed regions below the last unit-circle boundary are
due to Hu, Long, and Sun \cite{HLS}; the explicit arcs arise from Tschauner's
decoupling \cite{Tschauner,BKV}, and Wang and Zhou \cite{WZ} identified their
role as the boundaries of the complex-saddle region, with the spectrum
immediately below the lower arc changing type at the contact as shown below.

\begin{theorem}[Open regions and boundary multiplicities]\label{ii-sh-full-diagram}
On $0<e<1$, $0<\beta<9$, a nowhere-zero real-analytic function $H$
satisfies
\begin{equation}\label{ii-sh-full-disc}
 \Delta= g H^2,\qquad
 g=\frac{(\beta-b_-)(\beta-b_+)}{9-\beta}.
\end{equation}
In particular, $\Delta_\beta<0$ on $b_-$ and $\Delta_\beta>0$ on $b_+$.
All four multipliers in each of the following open regions are simple, and the
two elliptic pairs in the first and third rows give strong linear stability
under symplectic perturbations.
\begin{center}\small
\begin{tabular}{@{}lll@{}}
\toprule
Mass interval & Eccentricity & Essential spectrum\\
\midrule
$0<\beta<\beta_s$ & all & two elliptic pairs\\
$\beta_s<\beta<\beta_m$ & all & one elliptic and one negative pair\\
$\beta_m<\beta<b_-$ & $e<e_c$ & two elliptic pairs\\
$\beta_m<\beta<b_-$ & $e>e_c$ & two negative pairs\\
$b_-<\beta<b_+$ & all & a nonreal off-circle quartet\\
$b_+<\beta<9$ & all & two positive pairs\\
\bottomrule
\end{tabular}
\end{center}
Every complex eigenspace on the interior parameter domain has dimension one,
and the following table lists all boundary blocks there, with $J_k(z)$ denoting
a single Jordan block of size $k$.
\begin{center}\small
\begin{tabular}{@{}lll@{}}
\toprule
Boundary & Eccentricity & Jordan blocks or remaining simple pair\\
\midrule
$\beta_s$ & all & $J_2(-1)$ and one elliptic pair\\
$\beta_m$ & $e<e_c$ & $J_2(-1)$ and one elliptic pair\\
$\beta_m$ & $e>e_c$ & $J_2(-1)$ and one negative pair\\
$b_-$ & $e<e_c$ & $J_2(z),J_2(\bar z)$, $|z|=1$, $z\notin\R$\\
$b_-$ & $e>e_c$ & $J_2(z),J_2(z^{-1})$, $z<-1$\\
$b_+$ & all & $J_2(z),J_2(z^{-1})$, $z>1$\\
$b_-=\beta_m$ & $e=e_c$ & $J_4(-1)$\\
\bottomrule
\end{tabular}
\end{center}
\end{theorem}

\begin{proof}
For \(0<\beta<9\) the stiffness in \eqref{ii:model} is positive definite and
nonscalar; periodic nondegeneracy \cite{HLS} and the geometric simplicity of
all other multipliers \cite[Theorem~B]{HWmultiplicity} (recalled for this
family in Appendix~\ref{ii-app:simple}) imply that every eigenvalue has a
single Jordan block, so the minimal polynomial of \(M\) has degree four and no
identity \(M^2-hM+I=0\) can hold.

We next show that the factor $H$ in \eqref{ii-sh-full-disc} stays nonzero on
the arcs themselves. Write $P_c=P_0+cL$ in \eqref{ii-sh-ric}, where
$L=-J/(4q)$ is invertible, and for arbitrary complex $c$ let $h(c)$ be the
trace of the return of $u'=P_c u$; its determinant is one since
$\operatorname{tr}P_c=q'/q$, and $h$ is entire in $c$ and locally holomorphic
in the physical parameters. Whenever $c^2=g\ne0$, the two periodic graphs of
$P_c$ and $P_{-c}$ decompose the physical return with reduced traces $h(c)$ and
$h(-c)$, so
\[
 a_0=\tfrac12(h(c)+h(-c)),\qquad H=(h(c)-h(-c))/c
\]
are entire in $c^2$ and become single-valued real-analytic functions upon
substituting $c^2=g$, giving $\Delta=gH^2$ by continuity. In position--velocity
coordinates let $\widehat M$ be the return, and set $P=P_0(0)$ and $L=L(0)$ in
\[
 T=\begin{pmatrix}I&I\\P+cL&P-cL\end{pmatrix},\qquad
 Z=cT\begin{pmatrix}I&0\\0&-I\end{pmatrix}T^{-1}
  =\begin{pmatrix}-L^{-1}P&L^{-1}\\
       c^2L-PL^{-1}P&PL^{-1}\end{pmatrix}.
\]
The Cayley--Hamilton identity on each graph gives
$\widehat M^2-a_0\widehat M+I=(H/2)\widehat M Z$, and since the right-hand
formula for $Z$ extends through $c=0$, the identity persists there as well.
Because any zero of $H$ would produce the quadratic annihilator excluded above,
$H$ never vanishes (even at the fourfold contact), and differentiating $g$ yields
the two asserted derivative signs.

When $g<0$, the reduced traces are distinct complex conjugates and the spectrum
is a nonreal off-circle quartet; when $g>0$, both traces are real and distinct,
lying in $(-2,2)$ for an elliptic pair, below $-2$ for a negative pair, and
above $2$ for a positive pair. No trace can cross $2$, and a trace reaches $-2$
precisely along $\beta_s$ and $\beta_m$. On $b_-$ the common trace lies in
$(-2,2)$ before contact and below $-2$ after contact by the proof of
Theorem~\ref{ii-sh-contact}, whereas on $b_+$ it exceeds $2$; continuation
across each adjacent open real region yields the third, fourth, and last rows
of the first table, while the first two rows follow from \cite{HLS} (with the
real pair in the mixed region negative by continuation from $-1$). Because a
simple nonreal unit pair persists on the circle under small symplectic
perturbations (as its conjugate and reciprocal must coincide locally), the
first and third rows are strongly linearly stable. Along an explicit arc the
characteristic polynomial is $(z^2-hz+1)^2$, so the trace values above together
with geometric simplicity give two Jordan blocks of size two, or one block of
size four at contact; at every remaining antiperiodic point we have
$\Delta\ne0$, so $-1$ has algebraic multiplicity two and forms one size-two
block while the complementary pair is fixed by continuity from the adjacent
regions, completing the second table.
\end{proof}

\paragraph{The regular edges of the parameter rectangle.}
At $e=0$, substituting $u=\exp(\sigma t)u_0$ into \eqref{ii:model} shows that
the exponents solve $\sigma^4+\sigma^2+\beta/4=0$. For $0<\beta<1$, they are
$\pm i\omega_\pm$ with $\omega_\pm^2=(1\pm\sqrt{1-\beta})/2$, and since
$1<\omega_++\omega_-<\sqrt2$ and $0<\omega_+-\omega_-<1$, no frequency
coincidences occur except at $\beta=3/4$, where $\pm i/2$ yield a semisimple
double $-1$; all four multipliers are otherwise simple and elliptic. At
$\beta=1$, the two double exponents $\pm i/\sqrt2$ each have a one-dimensional
kernel in the exponential ansatz matrix (as its off-diagonal entry is nonzero),
producing two unit multipliers each with a size-two Jordan block, while for
$1<\beta<9$ the exponents have nonzero real part and imaginary parts of
absolute value $\sqrt{1+\sqrt\beta}/2\in(1/\sqrt2,1)$, giving a simple nonreal
off-circle quartet.

At $\beta=9$, the periodic rotation from the proof of
Theorem~\ref{ii-sh-contact} produces two identical scalar Hill equations whose
positive quadratic form excludes all unit twists; because the scalar
determinant-one return has trace greater than $2$ at $e=0$, it stays above $2$
for all $0\le e<1$, so its two positive reciprocal multipliers are each
semisimple of multiplicity two in the full essential return. At $\beta=0$,
$D_0=\operatorname{diag}(3,0)$ is the planar Kepler variational block: in
inertial variables, letting $\Phi_t$ be the Kepler flow with vector field $F$
and period $T=T(a)$, differentiating $\Phi_{T(a)}(X)=X$ across elliptic motions
gives $D\Phi_T=I-F\otimes dT$, whose nilpotent part has rank one and square
zero since $dT\ne0$ (by Kepler's third law) and $dT(F)=0$. Periodic rotation and
pulsation preserve this Jordan form (one $J_2(1)$ and two $J_1(1)$, including at
the circular orbit), while the edge $e=1$ is a singular limit of \eqref{ii:model}
whose two antiperiodic scales are determined in Section~\ref{ii-end-section}.

\section{The two scales at collision}\label{ii-end-section}

Put $d=1-e$, $h=\sqrt d$, $\rho=d^{1/4}$, and $k=(3-\sqrt{9-\beta})/2$, so
that $\beta=12k-4k^2$. We show that the two collision tangencies reported
numerically in \cite{MSS} admit convergent local parametrizations.

\begin{theorem}[Convergent power--log expansions at collision]\label{ii-end-theorem}
There exist holomorphic germs $T_A,T_B$ at the origin, real on real
arguments, such that the physical Lagrange roots satisfy
\begin{align*}
 k_A&=\rho T_A(\rho,\rho\log\rho),&
 T_A(0,0)&=\frac{1}{2^{1/4}\pi},\\
 k_B&=\rho^{10}T_B(\rho,\rho\log\rho),&
 T_B(0,0)&=\frac{2\sqrt2}{3\pi^2},
\end{align*}
and in fact $k_B=h^5\widetilde T_B(h,h\log h)$ for a holomorphic germ
$\widetilde T_B$. In the physical mass parameter,
\[
 \beta_A\sim\frac{12}{2^{1/4}\pi}d^{1/4},\qquad
 \beta_B\sim\frac{8\sqrt2}{\pi^2}d^{5/2}.
\]
In either expansion, truncating $T_j$ to total degree $N$ incurs an error of at
most $C Q^{N+1}/(1-Q)$, where $Q=\rho(1+|\log\rho|)/R<1$ and $C,R>0$ are
independent of $N$.
\end{theorem}

We use the forward angle $f\in[0,\pi]$ for the estimates and the reverse
angle $r=\pi-f$ for the connection. In the forward angle, sector $A$ carries
$y(0)=x'(0)=x(\pi)=y'(\pi)=0$ and sector $B$ carries
$x(0)=y'(0)=y(\pi)=x'(\pi)=0$. Throughout the $A$ calculation we retain a fixed
real parameter $\tau$ in the radial stiffness $3+(\tau-1)k$, with Lagrange
corresponding to $\tau=0$. Setting $p=x'$, $z=2x-y'$, and $q=1-e\cos r$ in the
reverse angle, the homogeneous shape equation takes the first-order form
\begin{equation}\label{ii-end-reverse}
 x'=p,\quad y'=2x-z,\quad
 p'=\big((3+(\tau-1)k)/q-4\big)x+2z,\quad z'=-ky/q .
\end{equation}
Denoting its transfer from $0$ to $\pi$ by $T$, the two sector determinants are
the principal minors with row and column sets $I_A=(y,p)$ and $I_B=(x,z)$,
respectively, since the free initial coordinates and the terminal coordinates
required to vanish are precisely $I_j$.

\subsection{Annular monodromy and the boundary determinants}

Set $s=\tan(r/2)$ and $\eta=\sqrt{(1-e)/(1+e)}$, and introduce the independent
variable $\chi=\eta/s$, so that differentiation with respect to $\log s$ along
an orbit of fixed $\eta$ is $\mathcal D=s\partial_s-\chi\partial_\chi$. We use
annular monodromy to isolate the leading boundary pairing.

\begin{lemma}[Holomorphic boundary determinants]\label{ii-end-transfer}
For $I=I_A,I_B$, the physical minor $D_I$ admits a representation
\begin{equation}\label{ii-end-minors}
 D_I=\eta^{-3}H_I(\eta,k,k\log\eta)
\end{equation}
with $H_I$ holomorphic near the origin in three independent variables.
At $\eta=0$, the upper minor factors into two holomorphic boundary
pairings, each divisible by $k$.
\end{lemma}

\begin{proof}
We first construct a frame at $k=0$. Let $a=1-e$, $b=1+e$, $\gamma=ab$,
$u=q\sin r$, and $v=uI$, where
\[
 I=\frac{s}{2b^2}-\frac{1}{2a^2s}
 -\frac{2e^3s}{\gamma^2(a+bs^2)}
 -\frac{6e^2}{\gamma^{5/2}}\big(\arctan(s/\eta)-\pi/2\big).
\]
Then $dI/dr=u^{-2}$ gives $uv'-u'v=1$. Setting $X=-q\cos r/e$,
$Y=\sin r\cos r-2\sin r/e$, and $V=(q^2I+\cot r)/e$, the columns
\begin{align*}
 E&=(0,1,0,0)^t,& G&=(u,q^2/e,u',0)^t,\\
 H&=(X,Y,X',1)^t,& W&=(v,V,v',0)^t
\end{align*}
solve \eqref{ii-end-reverse} at $k=0$ (for instance, $Y'=2X-1$, $V'=2v$, and
$u''=(3/q-4)u$), and their fundamental matrix $\Phi$ has determinant one.

Write $D_s=\operatorname{diag}(s^2,s,s,1)$,
$\Lambda=\operatorname{diag}(-1,1,0,-4)$, and $Q=D_s^{-1}\Phi s^{-\Lambda}$.
To verify that $Q$ is holomorphic at $(s,\chi)=(0,0)$, put $\eta=s\chi$ so
that $I=s^{-5}P$, where
\begin{align*}
 P&=\tfrac1{16}(1+s^2\chi^2)^2(2s^6+B_0/\chi^5),\\
 B_0&=-2\chi-\chi(1-s^2\chi^2)^3/(1+\chi^2)
 +3(1-s^2\chi^2)^2(1+s^2\chi^2)\arctan\chi .
\end{align*}
The numerator $B_0$ is odd in $\chi$ and its coefficients of $\chi$ and
$\chi^3$ vanish, so $P$ is holomorphic. Since
$q/s^2=2(1+\chi^2)/[(1+s^2\chi^2)(1+s^2)]$ is also holomorphic, the powers of
$s$ in $u,v,X,Y,V$ factor respectively as $s^3,s^{-2},s^2,s,s^{-1}$ times
holomorphic functions, and differentiation by $(1+s^2)\mathcal D/(2s)$ lowers
these powers by one. This establishes holomorphicity of $Q$ entry by entry,
with $\det Q=1$ and fourth row $e_3^t$; at $\chi=0$, the same formulas show
that every entry is even in $s$ once the powers are factored out.

Substituting the state $D_s Qg$ yields
\[
 \mathcal Dg=(\Lambda+kK)g,\qquad
 K=-\frac{1+s^2\chi^2}{1+\chi^2}Q^{-1}NQ,\qquad
 N_{31}=1-\tau,\quad N_{42}=1,
\]
with all other entries of $N$ zero. Choose $\delta>0$ inside the bidisk of
holomorphy of $K$, and let $T_g$ denote the transfer from $s=\eta/\delta$ to
$s=\delta$ and $M$ the positive monodromy based at $s=\delta$. On the fixed
outer circle the coefficient is holomorphic in $(\eta,k)$ up to $\eta=0$;
since $M(\eta,0)=\exp(2\pi i\Lambda)=I$, the local matrix logarithm defines a
holomorphic matrix $B=(2\pi i k)^{-1}\log M$.

We establish the holomorphy of the regularized transfer simultaneously for both
exterior representations $\mathcal R=\bigwedge^2$ and
$\bigwedge^2\circ\bigwedge^2$. Denoting the induced differential by
$\mathcal R_*$ and the largest diagonal entry of $\mathcal R_*\Lambda$ by $m$,
we claim that
\[
 \mathcal G_{\mathcal R}=\eta^m\exp\bigl((k\log\eta)\mathcal R_*B\bigr)
                         \mathcal R(T_g)
\]
extends holomorphically across $\eta=0$. Indeed, in $\xi=\log(s/\delta)$ the
fundamental matrix normalized at zero satisfies $Y(\xi+2\pi i)=Y(\xi)M$ by
periodicity, so analytic continuation around $\eta=0$ maps
$T_g=Y(\log(\eta/\delta^2))^{-1}$ to $M^{-1}T_g$ and keeps the displayed
product single-valued. On a principal branch, connect the inner and outer
circles by a radial segment and an arc of bounded angle: along the radius the
generator in log-radius is the real diagonal matrix $\mathcal R_*\Lambda$ plus
an $O(|k|)$ perturbation (uniformly in both cutoffs), so the norm differential
inequality yields $\|\mathcal R(T_g)\|\le C|\eta|^{-m-C|k|}$ while the arc
contributes a bounded factor and the exponential costs at most
$C|\eta|^{-C|k|}$. After shrinking the $k$-disk, $\mathcal G_{\mathcal R}=O(|\eta|^{-1/2})$
uniformly, so Cauchy's Laurent formula eliminates all negative powers and gives
joint holomorphy.

Put $\mathcal G=\mathcal G_{\wedge^2}$. Because the two largest weights on
$\bigwedge^2$ are $1$ (on $h=e_2\wedge e_3$) and $0$, both $\bigwedge^2$ and
$\bigwedge^2\circ\bigwedge^2$ have $m=1$, giving
$\bigwedge^2\mathcal G=\eta\,\mathcal G_{\wedge^2\circ\wedge^2}$ and hence
$\operatorname{rank}\mathcal G(0,k)\le 1$. At $k=0$ the transfer is
$(\delta^2/\eta)^\Lambda$, so $\mathcal G(0,0)=\delta^2 hh^t$; its $hh$ entry
remains nonzero nearby, yielding the holomorphic rank-one factorization
$\mathcal G(0,k)=c(k)r(k)$ from its $h$-column and normalized $h$-row.

On the outer cap $[2\arctan\delta,\pi]$ the system has a regular holomorphic
transfer $O$, whereas on the inner cap setting $s=\eta v$ and scaling the state
by $D_\eta$ yields
\begin{align*}
 \dot x_0&=2p_0/(1+\eta^2v^2),\qquad
 \dot y_0=(4\eta^2x_0-2z_0)/(1+\eta^2v^2),\\
 \dot p_0&=\frac{(3+(\tau-1)k)(1+\eta^2)}{1+v^2}x_0+
       \frac{4z_0-8\eta^2x_0}{1+\eta^2v^2},\\
 \dot z_0&=-k(1+\eta^2)y_0/(1+v^2),
\end{align*}
where dots denote $d/dv$, with a holomorphic Picard transfer $U$ on
$[0,1/\delta]$. Denoting the cut evaluations of $Q$ by $Q_i,Q_o$ and putting
$C_i=Q_i^{-1}D_\delta U$ and $C_o=O D_\delta Q_o$, the identity
$D_{\eta/\delta}^{-1}D_\eta=D_\delta$ expresses the physical transfer as
$C_o T_g C_i D_\eta^{-1}$. Writing $h_I=e_i\wedge e_j$ for $I=(i,j)$, both
physical input planes have total $D_\eta$-weight two, so
\[
 H_I=h_I^t\Bigl(\bigwedge^2 C_o\Bigr)\exp\bigl[-Z(\bigwedge^2)_*B\bigr]\,
                   \mathcal G\Bigl(\bigwedge^2 C_i\Bigr)h_I
\]
satisfies \eqref{ii-end-minors}, and at $\eta=0$ the factorization
$\mathcal G=cr$ separates its output and input pairings.

At $(\eta,k)=(0,0)$ the input row $r$ is proportional to $h^t$. Since the
fourth row of $Q_i$ is $e_3^t$ and $z_0$ is constant along the unperturbed
inner cap, we have $e_3^t C_i=e_4^t$; because the upper input plane
$I_A=(y,p)$ has zero fourth coordinate, its input pairing vanishes at $k=0$.

For the output pairing we must allow arbitrary independent $Z$. Differentiating
the outer-circle monodromy at $(\eta,k)=(0,0)$ gives
\[
 B(0,0)=\delta^\Lambda
 \operatorname*{Res}_{s=0}
 \left(s^{-\Lambda}K(s,0)s^\Lambda\,\frac{ds}{s}\right)
 \delta^{-\Lambda},
\]
which in column three extracts degrees $-1,1,0,-4$, respectively. Because the
odd coefficient vanishes by the parity of $Q(s,0)$, $B(0,0)$ preserves the line
$\mathbb C e_3$, so $\exp[-Z(\bigwedge^2)_*B(0,0)]h$ retains $e_3$ as a wedge
factor for every $Z$. At the unperturbed outer cut, $D_\delta Q_o e_3$ is the
column $H$ of the explicit frame, whose terminal value $(2,0,0,1)^t$ has
vanishing projection onto $(y,p)$. The output pairing at $k=0$ therefore
vanishes identically in $Z$, so each pairing carries a holomorphic factor of
$k$, proving the lemma.
\end{proof}

Taylor division in $\eta$ now gives
\begin{equation}\label{ii-end-div}
 H_{I_A}=\eta H_1(\eta,k,Z)+k^2 H_2(k,Z).
\end{equation}
At $k=0$, the free initial $y$ and $p$ columns are $E$ and $(G-a^2E/e)/a$,
whose terminal minor is $-b/a=-\eta^{-2}$, so in particular $H_1(0,0,0)=-1$.

\subsection{The central equation for the upper root}

The factorization \eqref{ii-end-div} leaves a single constant to determine,
which we compute from the Green kernel of the regular forward problem. In this
subsection and the next, $q(f)=1+e\cos f$, $b=1+e$, $q_*(f)=1+\cos f$, and
$J=\left(\begin{smallmatrix}0&-1\\1&0\end{smallmatrix}\right)$.

\begin{lemma}[Leading coefficient of the upper root]\label{ii-end-rank}
For fixed $\tau>-4$, the separated forward $A$ problem with stiffnesses
$(3+(\tau-1)k,k)$ admits a positive zero satisfying $k/\rho\to c_\tau$, where
\[
 c_\tau=\frac{2^{3/4}}{\pi\sqrt{\tau+4}} .
\]
\end{lemma}

\begin{proof}
Let $K_0f(t)=\int_0^t f(s)\,ds$, $V=q^{-1/2}$, and
$D=\operatorname{diag}(\tau-1,1)$, and work on $L^2((0,\pi);\mathbb C^2)$
with the $H^2$ domain defined by the $A$ boundary conditions. Write $R$ for
the inverse at $k=0$ of
$-\partial_f^2-2J\partial_f+q^{-1}\operatorname{diag}(3,0)$, and let $H$ be
the inverse of $-\partial_f^2+3/q-4$ under $x'(0)=x(\pi)=0$. The normalized
solutions of the homogeneous scalar equation are
\begin{equation}\label{ii-end-pair}
 u=q\sin f/b,\qquad
 g=b\left(\frac{\cos f}{q}-2e\sin f\,q\int_0^f\frac{\cos v}{q(v)^3}\,dv\right),
\end{equation}
with Cauchy data $(0,1)$ and $(1,0)$ at $f=0$, respectively, and terminal
values $g(\pi)=-b/d$, $u'(\pi)=-d/b$, so that
$H(t,s)=-g(\min(t,s))u(\max(t,s))$. Solving the second component of the
two-component system for source $(f_1,f_2)$ gives $y'+2x=K_0^*f_2$,
$x=H(f_1-2K_0^*f_2)$, and $y=K_0(K_0^*f_2-2x)$, whence
\[
 R=\begin{pmatrix} H&-2HK_0^*\\-2K_0H&K_0K_0^*+4K_0HK_0^*\end{pmatrix}.
\]

In the estimates below, all constants are uniform for $1/2\le e<1$. Writing
$x=\pi-t$ and $y=\pi-s$, we have $q(t)\asymp d+x^2$ and
$\int_0^t q^{-p}\le C_p q(t)^{1/2-p}$ for $p>1/2$. Since
$|\sin t|\le C\sqrt{q(t)}$, \eqref{ii-end-pair} gives $|g(t)|\le C/q(t)$, and
monotonicity of $q$ implies $|H(t,s)|\le C\min(x,y)$. Setting
$v(t)=R(t,\pi)e_2$ and $r_e=v_2(\pi)$ with $e_2=(0,1)^t$, we obtain
\begin{align*}
 v&=(-2H1,\ t+4K_0H1)^t,\qquad
 H1=(q-q^2)/(2e^2)+d^2g/(2eb),\\
 |r_e|&=\left|\frac{2d^2}{eb}\int_0^\pi g\right|\le C d^{3/2},
\end{align*}
because the scalar operator applied to $q-q^2$ gives $2e^2$ and
$\int_0^\pi(q-q^2)=-e^2\pi/2$. This formula also yields
$\|v-v_*\|_\infty\le Cd$ for
\[
 v_*=(q_*\cos t,-\sin t(2+\cos t))^t,\qquad \phi=q_*^{-1/2}v_* .
\]
Since $q\ge e q_*$ and $\phi$ is bounded, dominated convergence gives
$Vv_*\to\phi$ in $L^2$, while $\|V(v-v_*)\|_2\le Cd(\int q^{-1})^{1/2}\to 0$.

The unbounded part of $VRV$ comes from its boundary column and row. Subtracting
these traces via $E(t,s)=R(t,s)-v(t)e_2^t-e_2v(s)^t$ gives $E_{11}=H$,
$E_{12}=2\int_s^\pi H(t,z)\,dz$ (with $E_{21}(t,s)=E_{12}(s,t)$), and
\[
 E_{22}=\min(x,y)+4\int_t^\pi\int_s^\pi H(v,z)\,dz\,dv-r_e .
\]
Its nonconstant entries are bounded by $C\min(x,y)$, so after weighting by $V$
on both sides their kernels are bounded by $C/\max(x,y)$ (using $q(t)\ge c x^2$),
which defines a bounded operator on $L^2$ by Schur's test with weight
$y^{-1/2}$ (since the $y$-integral is at most $4x^{-1/2}$ and symmetrically in
$x$), while the constant term has weighted norm $|r_e|\int q^{-1}=O(d)$; thus
$\|VEV\|=O(1)$.

The single concentrating vector $Ve_2$ has squared norm
$\int_0^\pi q^{-1}=\pi/\sqrt{db}$, whereas $Vv\to\phi$, fixing the coupling
scale $\rho=d^{1/4}$ between the boundary trace and the outer solution. Putting
$a_e=\rho Ve_2$, its squared norm $\mu_e=\pi/\sqrt b$ stays bounded and
$a_e\rightharpoonup 0$ (as $a_e\to 0$ uniformly on compact subsets of
$[0,\pi)$), so $\alpha_e=\langle a_e,\phi\rangle\to 0$. With
$(a\otimes b)u=a\langle b,u\rangle$ and $Da_e=a_e$, the trace decomposition
gives in operator norm
\begin{equation}\label{ii-end-greenrank}
 S_e:=\rho D V R V=(D\phi)\otimes a_e+a_e\otimes\phi+o(1).
\end{equation}
One further scalar involves the entire outer solution:
\[
 \nu_\tau=\langle\phi,D\phi\rangle
 =\int_0^\pi\{4-4\cos^2t-2\cos^3t+
                 \tau(1+\cos t)\cos^2t\}\,dt=\pi(\tau+4)/2.
\]

Substituting $R$ and the boundedly invertible weight $V$ into the finite-$e$
equation shows that the problem at $k=\rho t$ has a kernel precisely when
$I+tS_e$ is singular. Writing $L_e$ for the $o(1)$ remainder in
\eqref{ii-end-greenrank}, $I+tL_e$ is invertible uniformly for large $e$ on
bounded sets of complex $t$, and factoring it out reduces the determinant to a
$2\times 2$ determinant that equals, up to $o(1)$ locally uniformly in $t$,
\[
 \det\begin{pmatrix}1+t\alpha_e&t\mu_e\\t\nu_\tau&1+t\alpha_e\end{pmatrix}
 \longrightarrow 1-\pi^2(\tau+4)t^2/(2\sqrt2),
\]
where the uniform error follows from boundedness of the two column vectors and
two functionals in \eqref{ii-end-greenrank}. Because the reduced determinant is
real for real $t$, on any sufficiently small closed interval around the simple
positive root $c_\tau$ of the limit its endpoint signs are opposite for $e$
close to one, giving a zero there by continuity whose convergence to $c_\tau$
follows from uniform convergence.
\end{proof}

Now substitute
\[
 \eta=\rho^2 g_0(\rho),\quad k=\rho t,\quad
 Z=2t\zeta+\rho t\log g_0(\rho),\qquad
 g_0(\rho)=(2-\rho^4)^{-1/2},
\]
with physical restriction $\zeta=\rho\log\rho$. By \eqref{ii-end-div}, the
scaled function $F_A=H_{I_A}/\eta=H_1+t^2 H_2/g_0$ is holomorphic even for
independent $(\rho,\zeta,t)$, with central value $-1+\sqrt2 H_2(0,0)t^2$.
Since the arguments tend to $(0,0,c_\tau)$ along the zeros of
Lemma~\ref{ii-end-rank}, we obtain
\begin{equation}\label{ii-end-centralA}
 F_A(0,0,t)=-1+t^2/c_\tau^2 .
\end{equation}
The holomorphic implicit function theorem therefore produces a unique germ
$t=T_\tau(\rho,\zeta)$ with $T_\tau(0,0)=c_\tau$, real on real arguments by
conjugation; for positive $\rho$ its restriction yields an actual separated
zero by \eqref{ii-end-minors} with nonzero derivative in the scaled parameter.
Setting $\tau=0$ proves the assertion for branch $A$ in
Theorem~\ref{ii-end-theorem}.

\subsection{The central equation for the lower root}

Take $\tau=0$ and put $H_B=H_{I_B}$ in \eqref{ii-end-minors}. The required
derivative can be computed directly in the forward angle without introducing a
Green operator for the singular endpoint problem.

\begin{lemma}[Leading jet of the lower determinant]\label{ii-end-Bjet}
The actual forward $B$ determinant $D_B$ satisfies
\[
 D_B(0,e)=-\eta^2,\qquad
 \lim_{e\uparrow1}\eta^3\partial_k D_B(0,e)=3\pi^2/16,
\]
as well as $D_B=\eta^{-3}H_B(\eta,k,k\log\eta)$.
\end{lemma}

\begin{proof}
In the state coordinates $(x,y,x_f',w)$ with $w=2x+y_f'$, let $F$ denote the
forward transfer, so that $D_B$ is its principal minor $(y,x_f')$. Reversal at
the two ends sends the state to $(x,y,-x_f',w)$, giving $T=S F^{-1}S$ with
$S=\operatorname{diag}(1,1,-1,1)$. Since the first-order equation is traceless,
$\det F=1$, and Jacobi's complementary minor identity yields the final
statement.

At $k=0$, the free initial columns have position components $(0,1)$ and $(u,v)$
with $u$ given by \eqref{ii-end-pair} and $v'=-2u$, giving
$D_B(0,e)=u'(\pi)=-d/b=-\eta^2$.

Fix $e<1$. Let $P$ extract $(y,x_f')$ from the state, let $C$ extract $(x,y)$,
and set $N(a,b)=(0,0,a,b)^t$. Denote the forward transfer at $k=0$ by $F(t,s)$
and put $S_0=PF(\pi,0)P^t$. With $V=q^{-1}\operatorname{diag}(-1,1)$, the
$k$-derivative of the state matrix is $NVC$, whereas a source $h$ on the right
of the second-order equation enters the state system as $-Nh$; for the inverse
under the $B$ conditions, variation of constants and the terminal condition
$PX(\pi)=0$ therefore yield the position kernel
\[
 G(t,s)=CF(t,0)P^t S_0^{-1}PF(\pi,s)N
            -\mathbf1_{s<t}CF(t,s)N .
\]
Because $CN=0$, the second term vanishes on the diagonal. Differentiating the
transfer by variation of constants and applying Jacobi's determinant formula
together with cyclicity of the finite-dimensional trace gives
\[
 \frac{\partial_kD_B(0,e)}{D_B(0,e)}
 =\int_0^\pi\operatorname{tr}\!\left(
 S_0^{-1}PF(\pi,t)NV(t)CF(t,0)P^t\right)dt
 =\int_0^\pi\operatorname{tr}(G(t,t)V(t))\,dt .
\]

To isolate the leading contribution in this trace, set
$Kf(t)=\int_0^t f(s)\,ds$. For a source $(h_1,h_2)$ the second equation gives
$w=-Kh_2$ and $y=K^*(2x+Kh_2)$ (since $w(0)=y(\pi)=0$), reducing the first
equation to $(-\partial_f^2+3/q-4)x=h_1+2Kh_2$ under $x(0)=x'(\pi)=0$. By
\eqref{ii-end-pair}, the scalar inverse has kernel
\[
 H(t,s)=u(\min(t,s))g(\max(t,s))
         -\frac{g'(\pi)}{u'(\pi)}u(t)u(s)
\]
(as the right solution has zero derivative at $\pi$ and the derivative jumps by
$-1$ across $t=s$), so the full inverse equals
$\binom{I}{2K^*}H\begin{pmatrix}I&2K\end{pmatrix}+\operatorname{diag}(0,K^*K)$,
with rank-one term generated by $(u,w_0)$ where
$w_0=2K^*u=(1+\cos f)(q+d)/b$.

Because $u$ is bounded and $|g|\le C/q$, the remaining scalar kernel is
$O(d^{-1})$ uniformly on the square, a bound preserved by integration on either
side, while $K^*K$ has a bounded kernel; since
$\int_0^\pi q^{-1}=O(d^{-1/2})$, the remainder of the diagonal trace integral
is $O(d^{-3/2})$. Using $D_B(0,e)=u'(\pi)$, we obtain
\[
 \eta^3\partial_kD_B(0,e)=
 -\eta^3 g'(\pi)\int_0^\pi (w_0^2-u^2)/q\,df+O(d).
\]
Differentiating $\int_0^\pi q^{-2}=\pi/(db)^{3/2}$ with respect to $e$ in
\eqref{ii-end-pair} gives $g'(\pi)=-3\pi e^2/(db)^{3/2}$, whence
$-\eta^3g'(\pi)\to 3\pi/8$. Because $u/q$ and $w_0/q$ are uniformly bounded,
dominated convergence evaluates the integral limit as
$\int_0^\pi \tfrac12\cos f(1+\cos f)^2\,df=\pi/2$, proving the lemma.
\end{proof}

By Lemma~\ref{ii-end-Bjet} and the identity theorem, $H_B(\eta,0,0)=-\eta^5$,
and finiteness of the physical derivative limit gives
\[
 H_{B,k}(\eta,0,0)+(\log\eta)H_{B,Z}(\eta,0,0)\longrightarrow
 a_0=3\pi^2/16 .
\]
Because both coefficient functions are analytic at the origin, this finiteness
forces $H_{B,Z}(0,0,0)=0$, $H_{B,Z}(\eta,0,0)=\eta b_1(\eta)$, and
$H_{B,k}(0,0,0)=a_0$. Taylor's formula in the independent variables $k,Z$
therefore yields
\[
 H_B=-\eta^5+k a_1(\eta)+Z\eta b_1(\eta)
      +k^2R_{20}+k Z R_{11}+Z^2R_{02},
\]
where $a_1(0)=a_0$ and the $R_{ij}$ are holomorphic. Substituting $k=\eta^5 t$
and $Z=\eta^4 t\xi$ makes every term divisible by $\eta^5$ in the convergent
power-series ring, with quotient satisfying
\begin{equation}\label{ii-end-centralB}
 F_B(0,0,t)=-1+(3\pi^2/16)t .
\end{equation}
The implicit function theorem gives a holomorphic germ $t=T(\eta,\xi)$ with
$T(0,0)=16/(3\pi^2)$, whose restriction to $\xi=\eta\log\eta$ for positive
$\eta$ yields a simple zero of the physical determinant.

To pass back to the common parameter $\rho$, set $g(h)=(2-h^2)^{-1/2}$, so
that $\eta=hg(h)$ and $\eta\log\eta=g(h)(h\log h)+hg(h)\log g(h)$, and define
\begin{align*}
 \widetilde T_B(h,v)&=g(h)^5 T\big(hg(h),g(h)v+hg(h)\log g(h)\big),\\
 T_B(\rho,\zeta)&=\widetilde T_B(\rho^2,2\rho\zeta).
\end{align*}
Since the two constructed positive zeros lie in different reflection sectors
and exhibit distinct leading orders, they correspond to the lower ($B$) and
upper ($A$) Lagrange branches. Finally, on a closed bidisk where $T_A$ and
$T_B$ are bounded, Cauchy's inequalities bound each bidegree-$(p,q)$ Taylor
coefficient by $CR^{-p-q}$; summing over $p+q=n$ via
$\sum_{p+q=n}|\rho|^p|\zeta|^q\le(|\rho|+|\zeta|)^n$ and then summing the
geometric tail yields the stated truncation bound, completing the proof of
Theorem~\ref{ii-end-theorem}.

\subsection{The alternating mode of a regular polygon}

The parameter $\tau$ in Lemma~\ref{ii-end-rank} also arises in a family of
Newtonian shape blocks. Fix an even integer $n\ge4$, place unit masses at the
vertices $Q_j=(\cos\theta_j,\sin\theta_j)$ ($\theta_j=2\pi j/n$) and a mass
$M>0$ at the center, set $\ell_j=2\sin(\pi j/n)$ for $1\le j<n$, and define
\[
 A_n=\sum_{j\ \mathrm{even}}\frac{1-\cos\theta_j}{\ell_j^3},\quad
 B_n=\sum_{j\ \mathrm{odd}}\frac{1+\cos\theta_j}{\ell_j^3},\quad
 C_n=\sum_{j\ \mathrm{odd}}\frac{1-\cos\theta_j}{\ell_j^3}.
\]

\begin{corollary}[Alternating mode of a regular polygon]\label{ii-end-polygon}
For fixed even $n\ge4$, set
\[
 k=\frac{C_n+2B_n}{M+A_n+C_n},\qquad
 \tau_n=\frac{B_n-C_n}{C_n+2B_n}.
\]
In the alternating block, the reflection sector $A$ admits a simple local curve
$k=\rho T_{\tau_n}(\rho,\rho\log\rho)$ with
$T_{\tau_n}(0,0)=2^{3/4}/(\pi\sqrt{\tau_n+4})$, whose central mass
$M=(C_n+2B_n)/k-A_n-C_n$ is positive for all sufficiently small $\rho>0$. Here
simplicity refers to this block.
\end{corollary}

\begin{proof}
Tangential forces cancel in pairs, giving the central-configuration multiplier
$\lambda=M+A_n+C_n$. Put $v_j=(-1)^j Q_j$ and $v_{\rm center}=0$; by the
finite geometric sums, both $v$ and $Jv$ are mass-orthogonal to translations,
$Q$, and $JQ$. The symmetry shifting labels by one and rotating by $-2\pi/n$
has real $-1$ eigenspace precisely $\operatorname{span}(v,Jv)$ (since its
$j$th component is $(-1)^j R_{\theta_j}v_0$ and the center must stay fixed),
and it commutes with both the mass Hessian and $J$, making this plane
invariant. Reflection across the horizontal axis combined with $j\mapsto-j$
fixes $v$ and negates $Jv$, so the selected reflection domain coincides with
the forward $A$ domain in the oriented frame $(v,Jv)/\sqrt n$.

For the positive Newtonian potential, an edge has Hessian
$(3\widehat d\widehat d^t-I)/|d|^3$; at $Q_0$, the differences $v_0-v_j$ are
parallel to the edge for even $j$ and perpendicular for odd $j$. Pairing $j$
with $n-j$ yields radial eigenvalue $2A_n-B_n$ and tangential eigenvalue
$-A_n+2B_n$ for the ring Hessian, while the center contributes $2M$ and $-M$,
respectively (with zero center row because $\sum_j v_j=\sum_j Jv_j=0$). Adding
the identity to the mass Hessian divided by $\lambda$ therefore yields the
shape stiffnesses $3+(\tau_n-1)k$ and $k$. Since $B_n,C_n>0$, we have
$-1<\tau_n<1/2$, and the block carries the Kepler coefficient
$1/(1+e\cos f)$, so the construction through \eqref{ii-end-centralA} applies
with $\tau=\tau_n$ (on a convergence neighborhood that may depend on $n$).
\end{proof}
\appendix
\section{Geometric simplicity in the open Lagrange domain}
\label{ii-app:simple}

For completeness we give the complex-time argument for the special case
of \cite[Theorem~B]{HWmultiplicity} used in the stability table.

\begin{lemma}[Geometric simplicity on the open Lagrange domain]\label{ii-app:geom}
For \(0<e<1\) and \(0<\beta<9\), every Floquet multiplier of the essential
return has geometric multiplicity one, and the return admits no quadratic
annihilator of the form \(M^2-hM+I=0\).
\end{lemma}
\begin{proof}
Rotate the stiffness as in \eqref{ii-sh-ric} to $D=\operatorname{diag}(d_1,d_2)$ with
$0<d_1=\kappa<3-\kappa=d_2$, and set $z=\exp(it)$, $\vartheta=z\partial_z$,
$x=(x_1,x_2)^T$, $u=(x_1,i x_2)^T$, $X=\left(\begin{smallmatrix}0&1\\1&0\end{smallmatrix}\right)$,
$E=\operatorname{diag}(1,-1)$, and $y=E(\vartheta x-Xx)$, so that
\[
 \vartheta\binom{y}{x}
 =\begin{pmatrix}-X&E(I-D/q)\\ E&X\end{pmatrix}\binom{y}{x},
 \qquad q=1+e(z+z^{-1})/2 .
\]
Equip the phase space with the bilinear symplectic form
$\omega((y,x),(Y,Z))=x^TY-y^TZ$. Inside the unit circle the only singularities
are $0$ and $r=(-1+\sqrt{1-e^2})/e$: at $r$ the residue $B_r$ has the single
nonzero upper-right block $C_r=\operatorname{diag}(-d_1,d_2)/\sqrt{1-e^2}$,
whereas at $0$ the residue $A_0$ is diagonalizable with eigenspaces
$V_-,V_0,V_+$ of dimensions $1,2,1$ at exponents $-2,0,2$.

To compute the monodromy, take analytic real symplectic gauges equal to $I$ at
each singularity. At $r$ the normal form has constant residue $B_r$ because
$kI-\operatorname{ad}_{B_r}$ is invertible at every positive order $k$, whereas
at $0$ the homological equations eliminate every nonresonant component by
Hamiltonian generators and leave $z^2R_2+z^4R_4$ with $[A_0,R_k]=kR_k$.
For completeness, the formal gauge relating two analytic coefficient matrices
satisfies a linear regular singular system whose vectorized coefficients obey
$(nI-L_0)v_n=\sum_{j=1}^n L_jv_{n-j}$ with $\|L_j\|\le C R^{-j}$. Since $L_0$
represents $\operatorname{ad}_{A_0}$, we have $\|L_0\|\ge\rho(L_0)=4$; for
$n>2\|L_0\|$ both resonant orders are passed and $\|(nI-L_0)^{-1}\|\le 2/n$,
so induction yields convergence on any smaller disk and therefore preserves the
formal symplectic identity. The nilpotent part $N=R_2+R_4$ maps $V_0$ into
$V_+$, annihilates $V_+$, and maps $V_-$ into $V_0+V_+$; for any $v\in V_-$ we
have $N^2v\in V_+$ and $\omega(v,N^2v)=-\omega(Nv,Nv)=0$, which forces $N^2=0$
by the nondegenerate pairing between the two one-dimensional extreme spaces.

Continue both local real fundamental matrices to a common base point in
$(r,0)$, where their singular factors are $(z-r)^{B_r}$ and
$(-z)^{A_0}\exp(N\log(-z))$, respectively. Because a positive loop around both
singularities is conjugate to the physical period map, $M$ is similar to
$M_*=(I+iA)(I+iB)$, where $A$ and $B$ are the real symplectic conjugates of
$2\pi N$ and $2\pi B_r$, respectively; both are Hamiltonian and square to zero.
Periodic nondegeneracy \cite{HLS} gives $\ker A\cap\ker B=0$, so $A$ and $B$
each have rank two with complementary Lagrangian images. Rank two forces
$N(V_0)=V_+$, and hence there exist real based vectors $w,\phi$ satisfying
$Aw=\phi$, $A\phi=0$, and $\omega(w,\phi)=0$ (taking the preimage in $V_0$
before transport and noting that $V_0$ is symplectically orthogonal to $V_+$).
Using the $+2$ eigenvector $(1,-1,1,1)^T$, normalize $\phi$ so that its
position component is $z^2(1,1)^T+O(z^3)$ near $0$.

To evaluate its pairing at the second pole, set $s=-z$, $f=-1/q(-s)>0$,
$p=s\partial_s x_1$, $v=x_2-x_1$, and $k=s\partial_s v$, so that with a dot
denoting $s\partial_s$ the system reads
\[
 \dot x_1=p,\quad \dot p=2p+2k+fd_1x_1,\quad
 \dot v=k,\quad \dot k=-2k+fd_2v+f(d_2-d_1)x_1 .
\]
Substituting $f=2s/e+O(s^2)$ gives $v=2(d_2-d_1)s^3/(15e)+O(s^4)$, so all four
variables are strictly positive near $s=0$. Because the differential system has
nonnegative off-diagonal entries, it preserves strict positivity on every
compact subinterval of $(0,|r|)$, and in particular $x_1$ and $v$ are
increasing. The pole normal form yields finite position limits at $r$ (its
logarithmic term has vanishing position component and the analytic gauge is
$I+O(z-r)$), and evaluating in that normal form before transporting back gives
\[
 \omega(\phi,B\phi)
 =\frac{2\pi}{\sqrt{1-e^2}}
       \bigl(d_2 x_2(r)^2-d_1 x_1(r)^2\bigr)>0 ,
 \qquad x_2(r)>x_1(r)>0 .
\]
Since $M_*+M_*^{-1}=2I-AB-BA$, any quadratic relation $M_*^2-hM_*+I=0$ with
$h\in\mathbb C$, applied to $\phi$ and paired with $w$, would give
$0=(2-h)\omega(w,\phi)=-\omega(\phi,B\phi)$, a contradiction. In coordinates
adapted to the complementary images, $A$ has sole upper-right block $A_1$ and
$B$ has sole lower-left block $B_1$, both invertible, and $A_1 B_1$ cannot be
scalar since $AB+BA$ is not scalar. For $z\ne 1$, eliminating the lower
coordinate of $(M_*-zI)v=0$ identifies the eigenspace with
$\ker(A_1 B_1+(z-1)^2 I/z)$, which is at most one-dimensional because $A_1 B_1$
is a nonscalar $2\times 2$ matrix, while the multiplier $1$ is already
excluded.
\end{proof}

\subsection*{Declaration of generative AI and AI-assisted technologies in the manuscript preparation process}
During the preparation of this work, the authors used generative AI tools to
assist with language editing and the verification of algebraic calculations.
After using these tools, the authors reviewed and edited the manuscript as
needed and take full responsibility for the content of the publication.

% Note: replace 'arXiv preprint (2026)' in \bibitem{HWmultiplicity} with \href{https://arxiv.org/abs/2611.XXXXX}{arXiv:2611.XXXXX} (2026) once Paper I is announced.

{\small\smallskip\noindent
Xijun Hu and Li Wu\\
School of Mathematics, Shandong University, Jinan 250100, P. R. China\\
\href{mailto:xjhu@sdu.edu.cn}{\texttt{xjhu@sdu.edu.cn}}\quad
\href{mailto:vvvli@sdu.edu.cn}{\texttt{vvvli@sdu.edu.cn}}\par}

\begin{thebibliography}{99}\small\raggedright

\bibitem{BKV}
B. Boldizs\'ar, T. Kov\'acs and J. Vany\'o, A new perturbative solution to the motion around triangular Lagrangian points in the elliptic restricted three-body problem, \emph{Celest. Mech. Dyn. Astron.} \textbf{133} (2021), article 23. \href{https://doi.org/10.1007/s10569-021-10018-8}{doi:10.1007/s10569-021-10018-8}.

\bibitem{DanR}
J. M. A. Danby, Stability of the triangular points in the elliptic restricted problem of three bodies, \emph{Astron. J.} \textbf{69} (1964), 165--172.

\bibitem{DanG}
J. M. A. Danby, The stability of the triangular Lagrangian points in the general problem of three bodies, \emph{Astron. J.} \textbf{69} (1964), 294--296.

\bibitem{HLS}
X. Hu, Y. Long and S. Sun, Linear stability of elliptic Lagrangian solutions of the planar three-body problem via index theory, \emph{Arch. Ration. Mech. Anal.} \textbf{213} (2014), 993--1045. \href{https://doi.org/10.1007/s00205-014-0749-6}{doi:10.1007/s00205-014-0749-6}.

\bibitem{HOregion}
X. Hu and Y. Ou, An estimation for the hyperbolic region of elliptic Lagrangian solutions in the planar three-body problem, \emph{Regul. Chaotic Dyn.} \textbf{18} (2013), 732--741. \href{https://doi.org/10.1134/S1560354713060129}{doi:10.1134/S1560354713060129}.

\bibitem{HOS}
X. Hu, Y. Ou and J. Sun, Quantitative linear stability analysis of elliptic relative equilibria in the planar $N$-body problem, \href{https://arxiv.org/abs/2509.09809}{arXiv:2509.09809} (2025).

\bibitem{HOWtrace}
X. Hu, Y. Ou and P. Wang, Trace formula for linear Hamiltonian systems with its applications to elliptic Lagrangian solutions, \emph{Arch. Ration. Mech. Anal.} \textbf{216} (2015), 313--357. \href{https://doi.org/10.1007/s00205-014-0810-5}{doi:10.1007/s00205-014-0810-5}.

\bibitem{HSindex}
X. Hu and S. Sun, Morse index and stability of elliptic Lagrangian solutions in the planar three-body problem, \emph{Adv. Math.} \textbf{223} (2010), 98--119. \href{https://doi.org/10.1016/j.aim.2009.07.017}{doi:10.1016/j.aim.2009.07.017}.

\bibitem{HWmultiplicity}
X. Hu and L. Wu, \emph{Floquet multiplicities and Krein signatures of elliptic relative equilibria}, arXiv preprint (2026).

\bibitem{HWZ}
X. Hu, L. Wu and A. Zhu, Strict separation of antiperiodic degeneracy curves for elliptic Euler and Lagrangian solutions, \href{https://arxiv.org/abs/2611.01719v1}{arXiv:2611.01719v1} (2026).

\bibitem{MSSearly}
R. Mart\'inez, A. Sam\`a and C. Sim\'o, Stability of homographic solutions of the planar three-body problem with homogeneous potentials, in \emph{Proceedings of Equadiff 2003}, World Scientific (2005), 1005--1010.

\bibitem{MSSdiagram}
R. Mart\'inez, A. Sam\`a and C. Sim\'o, Stability diagram for 4D linear periodic systems with applications to homographic solutions, \emph{J. Differential Equations} \textbf{226} (2006), 619--651. \href{https://doi.org/10.1016/j.jde.2006.01.014}{doi:10.1016/j.jde.2006.01.014}.

\bibitem{MSS}
R. Mart\'inez, A. Sam\`a and C. Sim\'o, Analysis of the stability of a family of singular-limit linear periodic systems in \(\mathbb R^4\). Applications, \emph{J. Differential Equations} \textbf{226} (2006), 652--686. \href{https://doi.org/10.1016/j.jde.2005.09.012}{doi:10.1016/j.jde.2005.09.012}.

\bibitem{Rob02}
G. E. Roberts, Linear stability of the elliptic Lagrangian triangle solutions in the three-body problem, \emph{J. Differential Equations} \textbf{182} (2002), 191--218.

\bibitem{Tschauner}
J. Tschauner, Die Bewegung in der N\"ahe der Dreieckspunkte des elliptischen eingeschr\"ankten Dreik\"orperproblems, \emph{Celestial Mechanics} \textbf{3} (1971), 189--196. \href{https://doi.org/10.1007/BF01228032}{doi:10.1007/BF01228032}.

\bibitem{WZ}
J. Wang and Q. Zhou, Precise analytic boundaries and positive symplectic paths in the stability diagrams of elliptic relative equilibria, \href{https://arxiv.org/abs/2611.21172v1}{arXiv:2611.21172v1} (2026).

\end{thebibliography}
\end{document}